\documentclass[10pt]{amsart}
\usepackage{amssymb,amsmath,amsthm}
\usepackage{mathrsfs,dsfont,a4wide}
\theoremstyle{plain}
\newtheorem{theorem}{Theorem}[section]
\newtheorem{lemma}[theorem]{Lemma}
\newtheorem{proposition}[theorem]{Proposition}

\newtheorem{remark}[theorem]{Remark}
\newtheorem{definition}[theorem]{Definition}
\theoremstyle{definition}
\theoremstyle{remark}
\numberwithin{equation}{section}

\newcommand{\as}{{\mathcal A}}
\newcommand{\Os}{{\mathcal O}}
\newcommand{\hs}{{\mathcal H}}
\newcommand{\ks}{{\mathcal K}}

\newcommand{\fs}{{\mathcal F}}
\newcommand{\leb}{{\mathcal L}}

\newcommand{\bs}{{\mathcal B}}

\newcommand{\Ps}{{\mathcal P}}

\newcommand{\Es}{{\mathcal E}}

\newcommand{\rs}{{\mathcal R}}

\newcommand{\R}{{\mathbb R}}
\newcommand{\N}{{\mathbb N}}

\newcommand{\Msym}{{\rm M}^{2\times 2}_{\rm sym}}

\newcommand{\Om}{\Omega}
\newcommand{\Omb}{\overline{\Omega}}

\newcommand{\weakst}{\stackrel{\ast}{\rightharpoonup}}

\newcommand{\weak}{\rightharpoonup}
\newcommand{\wlystar}{$\text{weakly}^*$\;}

\newcommand{\eps}{\varepsilon}

\newcommand{\Div}{{\rm div}\,}

\title
[A density result for Sobolev spaces in dimension two]
{A density result for Sobolev spaces in dimension two, and applications to stability of nonlinear Neumann problems}
\author[A. Giacomini]
{Alessandro Giacomini}
\address[Alessandro Giacomini]{Dipartimento di Matematica, Facolt\`a di Ingegneria, Universit\`a degli Studi di Brescia, Via Valotti 9, 25133 Brescia, Italy}
\email[A. Giacomini]{alessandro.giacomini@ing.unibs.it}
\author[P. Trebeschi]
{Paola Trebeschi}
\address[Paola Trebeschi]{Dipartimento di Matematica, Facolt\`a di Ingegneria, Universit\`a degli Studi di Brescia, Via Valotti 9, 25133 Brescia, Italy}
\email[P. Trebeschi]{paola.trebeschi@ing.unibs.it}
\begin{document}
\begin{abstract}
We prove that if $\Om \subseteq \R^2$ is bounded and $\R^2 \setminus \Om$ satisfies suitable structural assumptions (for example it has a countable number of connected components), then $W^{1,2}(\Om)$ is dense in $W^{1,p}(\Om)$ for every $1\le p<2$. The main application of this density result is the study of
stability under boundary variations for nonlinear Neumann problems of the form
$$
\begin{cases}
-{\rm div}\,A(x,\nabla u)+B(x,u)=0      & \text{in }\Om, \\
A(x,\nabla u)\cdot \nu=0      & \text{on }\partial \Om,
\end{cases}
$$
where $A:\R^2\times \R^2 \to \R^2$ and $B:\R^2 \times \R \to \R$ are Carath\'eodory functions which satisfy standard monotonicity and growth conditions of order $p$.
\end{abstract}
\maketitle
\tableofcontents

\section{Introduction}
\label{intro}
In this paper we prove a density result for Sobolev spaces defined on two dimensional open bounded sets. More precisely, for $1 \le p<2$ and $\Om \subseteq \R^2$ open, bounded and belonging to the class $\as_p(\R^2)$ of admissible domains (see Definition \ref{maindef}), we prove that the Sobolev space $W^{1,2}(\Om)$ is dense in $W^{1,p}(\Om)$. The class $\as_p(\R^2)$ contains for example domains whose complements have a countable number of connected components or even whose complements are Cantor sets with small dimension.
\par
In the case $\Om$ is sufficiently regular (for example if it satisfies a cone condition), this density result is trivial because by means of extension operators and convolutions one can prove that $C^\infty(\Omb)$ is dense in $W^{1,p}(\Om)$. The situation is different when $\Om$ is irregular: extension operators cannot be employed, and the density of $C^\infty(\Omb)$ in $W^{1,p}(\Om)$ can fail, as in the case the domain contains a crack.
Even the density of $C^\infty(\Om)$ in $W^{1,p}(\Om)$ proved by Meyers and Serrin \cite{MeSe} which holds for every open bounded set $\Om$ cannot be used because  the control on the energy of order $2$ is available only well inside, and can be lost approaching the boundary.
In this direction, we refer the reader to the paper of O'Farrel \cite{Fa} for a counterexample to the density of $W^{1,\infty}(\Om)$ in $W^{1,p}(\Om)$ in the case $\Om$ is too irregular.
\par
The main motivation of our density result is the study of stability under boundary variations for two dimensional nonlinear Neumann problems of the form
\begin{equation}
\label{intrpbgeneral}
\begin{cases}
-{\rm div }A(x,\nabla u)+B(x,u)=0      &\text{in }\Om\\
A(x,\nabla u)\cdot \nu=0      &\text{on }\partial \Om,
\end{cases}
\end{equation}
where $A:\R^2\times \R^2 \to \R^2$ and $B:\R^2\times \R \to \R$ are Carath\'eodory functions satisfying standard monotonicity and growth conditions of order $p$
(see conditions \eqref{monotone}-\eqref{below}). Namely we are interested in the continuity of the map $\Om \to u_\Om$, where $u_\Om \in W^{1,p}(\Om)$ is the solution of \eqref{intrpbgeneral} in $\Om$ (see Section \ref{applicsec} for the precise sense of the continuity of this mapping). 
\par
The density of $W^{1,2}$ in $W^{1,p}$ is a key point to infer stability for problem \eqref{intrpbgeneral} from that of the linear equation
\begin{equation}
\label{intrpblaplacian2}
\begin{cases}
-\Delta u+u=f      &\text{in }\Om\\
\frac{\partial u}{\partial \nu}=0      &\text{on }\partial \Om.
\end{cases}
\end{equation}
Stability results for problem \eqref{intrpblaplacian2} have been obtained by several authors (see for example
\cite{BZ1}, \cite{BZ2}, \cite{BZ3},
\cite{ChD}, \cite{Che}, \cite{Buc}, \cite{Buc2}). 
These results hold in generic dimension $N$ under quite restrictive assumptions on $\Om$ and its possible perturbations. For example Chenais \cite{Che} proved stability for \eqref{intrpblaplacian2} under a uniform cone condition for the perturbed domains, and this condition excludes several interesting cases like those of domains containing cracks which are of interest in fracture mechanics. Moreover, the cone condition implies the existence of extension operators, and the density of $W^{1,2}$ in $W^{1,p}$ is trivial, so that the stability of \eqref{intrpbgeneral} holds under the same assumptions.
\par
In dimension $N=2$ the situation is different, and restrictions only on the topological nature of the domains have been individuated in order to achieve stability for \eqref{intrpblaplacian2}: this is the reason why we are interested in density for Sobolev spaces defined on two dimensional, possibly irregular, domains. Bucur and Varchon \cite{Buc} consider domains whose complements have a uniformly bounded number of connected components and prove that, if $\Om_n \to \Om$ in the Hausdorff complementary topology (see Section \ref{notations} for a definition), we have the stability $u_{\Om_n} \to u_\Om$ if and only if
$$
{\rm meas}(\Om_n) \to {\rm meas}(\Om).
$$
Under strict monotonicity assumptions for $A$ and $B$, Dal Maso, Ebobisse and Ponsiglione \cite{DMEP} proved
that the same conclusion holds for problem \eqref{intrpbgeneral} in the case $1<p<2$, while for $p>2$ stability is in general false (see \cite[Remark 3.7]{DMEP}).
The main tool they employ is the Mosco convergence of $W^{1,p}(\Om_n)$ to $W^{1,p}(\Om)$ (see Section \ref{notations} for a definition) which is equivalent to the stability of \eqref{intrpbgeneral} for every admissible $A$ and $B$. The Mosco convergence in the case $p=2$ is indeed a corollary of the stability result by Bucur and Varchon \cite{Buc}. Since they make use of conformal mappings, and these are not useful in a nonlinear setting, Dal Maso, Ebobisse and Ponsiglione provide a different proof of the Mosco convergence based on {\it nonlinear harmonic conjugates}. 
In view of our density result, the Mosco convergence when $1<p<2$ (and hence the stability result for \eqref{intrpbgeneral}) can be deduced from the case for $p=2$ (see Proposition \ref{moscoconv}).
\par 
In Section \ref{applicsec} we consider the nonlinear Neumann problems
\begin{equation}
\label{intrpb}
\begin{cases}
-{\rm div }A(x,\nabla u)+b(x)|u|^{p-2}u=h      &\text{in }\Om\\
A(x,\nabla u)\cdot \nu=0      &\text{on }\partial \Om,
\end{cases}
\end{equation}
where $b \in L^\infty(\R^2)$ is such that $b \ge 0$, and $h$ satisfies suitable assumptions in order to guarantee the existence of a solution.
These problems introduce some degeneracy with respect to problems \eqref{intrpbgeneral} because $b$ can vanish on subsets of $\Om$ with positive measure. As a consequence stability cannot be studied in terms of Mosco convergence of suitable functional spaces, because the two notions are in general not equivalent (see \cite[Remark 5.2]{Buc2}), and so in order to prove stability for \eqref{intrpb}, the results of Dal Maso, Ebobisse and Ponsiglione cannot be directly used.
\par
In the case $p=2$, and with $A(x,\xi)=\xi$, Bucur and Varchon \cite{Buc2} proved that if the complement of $\Om_n$ has a uniformly bounded number of connected components and $\Om_n \to \Om$ in the Hausdorff complementary topology, then 
stability holds if and only if
$$
{\rm meas}(\Om_n \cap \{b>0\}) \to {\rm meas}(\Om \cap \{b>0\}).
$$
We prove (Proposition \ref{stabilitynonlinear}) that the same result holds in the nonlinear case $1<p<2$. In the case $p>2$, stability does not hold in general (see \cite[Remark 3.7]{DMEP}).
\par
A second application of our density result is to a shape optimization problem, namely the optimal cutting of a membrane. The admissible cuts we consider are compact and connected sets which contain two given points. The case of a quadratic energy has been treated by Bucur, Buttazzo and Varchon in \cite{BBV}. In Proposition \ref{optimalcutprop} we prove the existence of an optimal cut for a nonlinear energy density $f(x,\xi)$ with growth of order $1<p \le 2$ in $\xi$. Moreover we prove a stability result for the associated Euler-Lagrange equation, which is of Neumann-Dirichlet type. We remark that in order to establish the existence of the optimal cut and the stability for the associated equation, the approximation results of Dal Maso, Ebobisse and Ponsiglione \cite{DMEP} in terms of Mosco convergence cannot be used (see Remark \ref{remcut}).
\par
Finally, in the Appendix, we show how the arguments of Section \ref{secgradients} provide a new proof of a result due to Chambolle \cite{Ch} concerning the density of $W^{1,2}$ in the space $LD^{1,2}$ of two dimensional linearized elasticity.
Our approach also covers the nonlinear case $LD^{1,p}$ for $1<p<2$. 
\par
The main step in the proof of our main result is given by Theorem \ref{mainthml1p}, which states the density of $W^{1,2}(\Om)$ in $W^{1,p}(\Om)$ at the level of the gradients. More precisely we prove that for every $u \in W^{1,p}(\Om)$ we have $\nabla u \in \overline H$, where
$$
H:=\{\nabla v\,:\, v \in W^{1,2}(\Om)\} \subseteq L^p(\Om,\R^2).
$$
We use the fact that $\overline H=(H^\perp)^\perp$, where $(\cdot)^\perp$ denotes the orthogonal in the sense of Banach spaces. Using Helmholtz Decomposition Theorem, in Lemma \ref{helmoltz} we characterize $H^\perp$ as the family of fields $\psi$ such that $R\psi=\nabla \phi$ with $\phi \in W^{1,p'}(\R^2)$ constant on the connected components of $\R^2 \setminus \Om$, where $p'$ is the conjugate exponent of $p$ and $R$ denotes a rotation of 90 degrees counterclockwise. Moreover, using the approximation given in Lemma \ref{image} and the fact that $\Om \in \as_p(\R^2)$, we can approximate $\phi$ through functions $\phi_n \in W^{1,p'}(\R^2)$ which are constant on a neighborhood of $\R^2 \setminus \Om$. Then the orthogonality of $\nabla u$ and $\psi$ follows by integration by parts.
\par
The paper is organized as follows. In Section \ref{notations} we introduce the notation and recall some useful notions employed in rest of the paper. Section \ref{secgradients} contains the density result (Theorem \ref{mainthmW1p}), while Section \ref{applicsec} contains the applications to stability of nonlinear Neumann problems and to the optimal cutting of a membrane.
In the Appendix we prove the density of $W^{1,2}$ in the spaces of planar elasticity.

\section{Notation and Preliminaries}
\label{notations}
In this section we introduce the basic notation and recall some notions employed in the rest of the paper.
\par
If $A \subseteq \R^N$ is open and $1 \le p \le +\infty$, we denote by $L^p(A)$ the usual space of $p$-summable functions on $A$ with norm indicated by $\|\cdot\|_p$. $W^{1,p}(A)$ will denote the Sobolev  space of functions in $L^p(A)$ whose gradient in the sense of distributions belongs to $L^p(A,\R^N)$, and we denote by $W^{1,p}_0(A)$ the closure in $W^{1,p}(A)$ of smooth functions with compact support.
\par
If $E \subseteq \R^N$, we will denote with ${\rm meas}(E)$ its $N$-dimensional Lebesgue measure, and by $\hs^\alpha(E)$ its $\alpha$-dimensional Hausdorff measure (see \cite[Chapter 2]{EG} for a definition). Moreover, we denote by $E^c$ the complementary set of $E$, and by $1_E$ its characteristic function, i.e., $1_E(x)=1$ if $x \in E$, $1_E(x)=0$ otherwise.

\vskip15pt\noindent 
{\bf Capacity.} Let $1<p<+\infty$, and let $E \subseteq \R^N$. We set
$$
c_p(E):=\inf\left\{\int_{\R^2} |\nabla u|^p+|u|^p\,dx\,:\,u \in W^{1,p}(\R^2), u \ge 1 \text{ a.e. on }E\right\}.
$$
For the properties of capacity, and its relevance in the theory of Sobolev spaces, we refer the reader to \cite{EG}. 
\par
We say that a property $\Ps(x)$ holds $c_p$-quasi everywhere (abbreviated $c_p$-q.e.) on a set $E \subseteq \R^N$ if it holds for every $x \in E$ except a subset $N$ of $E$ such that $c_p(N)=0$. 
\par
If $A \subseteq \R^N$ is open, every function $u \in W^{1,p}(A)$ admits a {\it quasicontinuous} representative, i.e., a representative $\tilde u$ such that for every $\eps>0$ there exists an open set $B_\eps$ with $c_p(B_\eps)<\eps$ and $\tilde u_{|A \setminus B_\eps}$ is continuous. Throughout the paper, we will identify a Sobolev function with its quasicontinuous representative. Notice that for $p>N$, the continuous representative of $u$ (which exists by Sobolev Embedding Theorem) is precisely the quasicontinuous representative. We will use the following fact: if $u_n \to u$ strongly in $W^{1,p}(A)$, we have that up to a subsequence $u_n \to u$ $c_p$-q.e. on $A$.
\par
The following lemma will be useful in Section \ref{secgradients} and in the Appendix
(a different proof can be obtained using the arguments contained in
\cite[Lemma 5.1, Lemma 5.2]{BT}).

\begin{lemma}
\label{continuous}
Let $u \in C(\R^2)$, $K \subseteq \R^2$ connected, and let $c \in \R$. If $u(x)=c$ for $c_2$-q.e. $x \in K$, then $u(x)=c$ for every $x \in K$.
\end{lemma}

\begin{proof}
By assumption we have that there exists $N \subseteq K$ such that $c_2(N)=0$ and
$u(x)=c$ for every $x \in K \setminus N$. If for every $x \in N$ there exists $x_n \in K\setminus N$ such that $x_n \to x$, by continuity of $u$ we conclude that also $u(x)=c$ and the result follows.
\par
By contradiction, let us assume that there exists $x \in N$ such that $x \not\in \overline{K \setminus N}$. Then there exists $\bar r>0$ such that $B(x,r) \cap (K \setminus N)=\emptyset$ for $r<\bar r$. Since $c_2(N)=0$, by \cite[Section 4.7.2, Theorem 4]{EG} 
we have that $\hs^\alpha(N)=0$ for every $\alpha>0$, and in particular $\hs^1(N)=0$. As a consequence, for every $0<\eps<\bar r$ we can find a covering $\{B(x_i,r_i)\}_{i \in \N}$ of $N$ such that $\sum_{i \in \N}r_i< \eps$. Let $\bs:=\cup_iB(x_i,r_i)$ and
$$
S:=\{r \in ]0,\bar r[\,:\, \partial B(x,r) \cap \bs\not=\emptyset\}.
$$
We have that ${\rm meas}(S)<\eps$, so that we can find $r<\bar r$ with $\partial B(x,r) \cap N=\emptyset$. Moreover, up to reducing $\eps$, we can assume that $N \setminus B(x,r) \not=\emptyset$, because otherwise we would get that $N=\{x\}$ with $x$ isolated from the rest of $K$, against its connectedness.
\par
Let us consider
$$
K_1:=K \cap \overline{B(x,r)}
\quad\text{and}\quad
K_2:=K \setminus B(x,r).
$$
$K_1$ and $K_2$ are closed in the relative topology of $K$. By construction they are not empty, disjoint and such that $K=K_1 \cup K_2$. But this is against the fact that $K$ is connected, and the proof is concluded. 
\end{proof}

\vskip20pt\noindent
{\bf Hausdorff metric on compact sets and Hausdorff complementary topology.}
Let $A$ be open and bounded in $\R^N$. We indicate the family of all compact subsets of $\overline A$ by $\ks(\overline A)$. $\ks(\overline A)$ can be endowed with the
Hausdorff metric $d_H$ defined by
$$
d_H(K_1,K_2) := \max \left\{ \sup_{x \in K_1} {\rm dist}(x,K_2), \sup_{y \in
K_2} {\rm dist}(y,K_1)\right\}
$$
with the conventions ${\rm dist}(x, \emptyset)= {\rm diam}(A)$ and $\sup
\emptyset=0$, so that $d_H(\emptyset, K)=0$ if $K=\emptyset$ and
$d_H(\emptyset,K)={\rm diam}(A)$ if $K \not=\emptyset$. It turns out that
$\ks(\overline A)$ endowed with the Hausdorff metric is a compact space
(see e.g. \cite{Ro}).
\par
In order to treat the stability of Neumann problems under boundary variations (see Section \ref{applicsec}), we will use the {\it Hausdorff complementary topology} on the family of open sets which is defined as follows.
Let $(\Om_n)_{n \in \N}$ be a sequence of open sets in $\R^N$. We say that $\Om_n \to \Om$ in the Hausdorff complementary topology if for every closed ball $B \subseteq \R^N$ we have 
$$
B \cap \Om_n^c \to B \cap \Om^c 
\quad\text{in the Hausdorff metric}.
$$

\vskip20pt\noindent
{\bf The Mosco convergence of Sobolev spaces.}
In Section \ref{applicsec}, we will refer to the notion of Mosco convergence of Sobolev spaces in connection with stability results for nonlinear Neumann problems. For the reader's convenience, we recall here the definition.
\par
Let $(\Om_n)_{n \in \N}$ be a sequence of uniformly bounded open subsets of $\R^N$, and let $1<p<+\infty$.
For every $u_n \in W^{1,p}(\Om_n)$, let us denote by $u_n1_{\Om_n}$ and by $\nabla u_n1_{\Om_n}$ the extension to zero outside $\Om_n$ of $u_n$ and $\nabla u_n$ respectively.
\par
If $\Om$ is a bounded open set in $\R^N$, we say that $W^{1,p}(\Om_n)$ converges to $W^{1,p}(\Om)$ in the sense of Mosco if the following two conditions hold.
\begin{itemize}
\item[$(M1)$] {\it Mosco-limsup condition.} For every $u \in W^{1,p}(\Om)$ there exists
$u_n \in W^{1,p}(\Om_n)$  such that
$$
\nabla u_n1_{\Om_n} \to \nabla u1_{\Om}
\quad\text{strongly in }L^p(\R^N,\R^N)
$$
and
$$
u_n1_{\Om_n} \to u1_{\Om}
\quad\text{strongly in }L^p(\R^N).
$$
\item[$(M2)$]{\it Mosco-liminf condition.} If $n_k$ is a sequence of indices converging to $+\infty$, $(u_k)_{k \in \N}$ is a sequence such that $u_k \in W^{1,p}(\Om_{n_k})$ for every $k$, and $u_k1_{\Om_{n_k}}$ converges weakly in $L^p(\R^N)$ to a function $\varphi$, while $\nabla u_k1_{\Om_{n_k}}$ converges weakly in $L^p(\R^N,\R^N)$ to a function $\Phi$, then there exists $u \in W^{1,p}(\Om)$ such that $\varphi=u1_\Om$ and $\Phi=\nabla u1_\Om$.
\end{itemize}
Using a diagonal argument, we have that in order to establish $(M1)$, it suffices to approximate functions belonging to a dense subset of $W^{1,p}(\Om)$.
This fact will be used several times in Section \ref{applicsec}.

\section{The density result}
\label{secgradients}
This section is devoted to the proof of the density of $W^{1,2}$ into $W^{1,p}$ with $1 \le p<2$ on a two-dimensional domain which satisfies a suitable structural assumption, for example if its complement has a countable number of connected components. Recall that the two-dimensional domain is not assumed to be regular (for example it may contain a crack), so that extension operators cannot be used. 
\par
First of all, we establish the density result at the level of the gradients (Theorem \ref{mainthml1p}). The extension to the full result on Sobolev spaces (Theorem \ref{mainthmW1p}) is then obtained through a truncation argument. 
\par
The class of admissible domains we consider is given in the following definition.

\begin{definition}{\bf (The class $\as_p(\R^2)$ of admissible domains)}
\label{maindef}
Let $1\le p<2$, and let $\Om \subseteq \R^2$ be open and bounded. Let $\{K_i\}_{i \in I}$ be the family of the connected components of $\Om^c$. We say that $\Om$ belongs to the class $\as_p(\R^2)$ of  admissible domains if for every $i \in I$ there exists $x_i \in K_i$ such that setting $E:=\{x_i, i \in I\}$ we have
\begin{equation}
\label{defOmadmis}
\hs^{2-p}(E)=0.
\end{equation}
\end{definition}

Notice that the class $\Os_l(\R^2)$ of two-dimensional domains such that their complements have at most $l$ connected components (which is relevant for stability of nonlinear Neumann problems, see Section \ref{stability}) is contained in
$\as_p(\R^2)$. Moreover $\as_p(\R^2)$ contains domains $\Om$ such that $\Om^c$ has a countable number of connected components, or even an uncountable number provided that there exists a suitable selection $E$ of $\{K_i\}_{i \in I}$ with zero Hausdorff measure of order $2-p$. We remark that condition \eqref{defOmadmis} is not referred to the connected components $K_i$ of $\Om^c$ but to a selection $E$ of $\{K_i\}_{i \in I}$: in particular it can be ${\rm meas}(K_i)>0$ (not only for the unbounded connected component).
\par
The following lemmas will be useful in the proof of Theorem \ref{mainthml1p}.

\begin{lemma}
\label{capzero}
Let $A \subseteq \R^2$ be open, and let $u \in W^{1,q}(A)$ with $q>2$. Then we have ${\rm meas}(u(E))=0$ for every $E \subseteq A$ such that $\hs^{\frac{q-2}{q-1}}(E)=0$
(in the case $q=+\infty$ we mean $\hs^1(E)=0$).
\end{lemma}

\begin{proof}
If $q=+\infty$, the result follows because $u$ is a locally Lipschitz function and
${\rm meas}(f(C))=\hs^1(f(C)) \le L \hs^1(C)$ for every $L$-Lipschitz function $f$ and every set $C$ (see \cite[Theorem 1, Section 2.4.1]{EG}).
\par
In the case $2<q<+\infty$, we follow the approach that Marcus and Mizel \cite{MM} developed to deal with $N$-property of Sobolev transformations (see \cite{FG} for a description of the problem of $N$-property, and \cite[Theorem 5.28]{FG}). 
\par
By Sobolev Embedding Theorem $u$ is a H\"older continuous function. Moreover, for any square $Q_r \subseteq A$ of side $r$ we have 
\begin{equation}
\label{sobolev}
|u(x)-\bar u_{Q_r}| \le C_q \|\nabla u\|_{L^q(Q_r,\R^2)}r^{1-2/q},
\end{equation}
where $\bar u_{Q_r}$ denotes the average of $u$ on $Q_r$, 
and $C_q$ depends only on $q$. From \eqref{sobolev} we deduce that $u(Q_r)$ is contained in an interval $I_{Q_r}$ of length at most
$$
l_{Q_r}:=2C_q \|\nabla u\|_{L^q(Q_r,\R^2)}r^{1-2/q}.
$$
\par
Let $E \subseteq A$ be such that $\hs^{\frac{q-2}{q-1}}(E)=0$, and let us fix $\eps>0$ and $\delta>0$. Since $\hs^{\frac{q-2}{q-1}}(E)=0$,
we can find a covering $\fs=\{Q_{r_i}(x_i)\}_{i \in \N}$ of $E$ with $Q_{r_i}(x_i) \subseteq A$,
\begin{equation}
\label{zeroalpha}
\sum_{i=0}^{+\infty}r_i^{\frac{q-2}{q-1}}<\eps
\end{equation}
and such that 
$$
l_i:=2C_q \|\nabla u\|_{L^q(Q_{r_i}(x_i),\R^2)}r_i^{1-2/q}<\delta.
$$
By Besicovich Covering Theorem (see \cite[Section 1.5.2, Theorem 2]{EG}) there exist
$m$ families $\fs_1, \fs_2,\dots,\fs_j,\dots,\fs_m \subseteq \fs$ of disjoint squares $\{Q_{r_{i,j}}(x_{i,j})\}_{i \in \N}$
such that 
$$
E \subseteq \bigcup_{j=1}^m \bigcup_{i=0}^{+\infty} Q_{r_{i,j}}(x_{i,j}).
$$
By H\"older inequality and by \eqref{zeroalpha} we deduce that
\begin{align*}
\sum_{j=1}^{m}\sum_{i=0}^{+\infty}l_{i,j}&=2C_q\sum_{j=1}^{m}\sum_{i=0}^{+\infty} \|\nabla u\|_{L^q(Q_{r_{i,j}}(x_{i,j}),\R^2)}r_{i,j}^{1-2/q} \\
&\le 2C_q \left( \sum_{j=1}^{m}\sum_{i=0}^{+\infty} \|\nabla u\|^q_{L^q(Q_{r_{i,j}}(x_{i,j}),\R^2)}\right)^{\frac{1}{q}} \left(\sum_{j=1}^{m}\sum_{i=0}^{+\infty}r_{i,j}^{\frac{q-2}{q-1}} \right)^{\frac{q-1}{q}} \\
&\le 2C_q m\|\nabla u\|_{L^q(A,\R^2)} \eps^{\frac{q-1}{q}}
\end{align*}
so that
$$
\hs^1_\delta(u(E)) \le 2C_q m\|\nabla u\|_{L^q(A,\R^2)} \eps^{\frac{q-1}{q}},
$$
where $\hs^1_\delta(E)$ denotes the $(1,\delta)$-Hausdorff pre-measure. Since $\hs^1(u(E))=\lim_{\delta \to 0}\hs^1_\delta(u(E))$, and $\hs^1(u(E))={\rm meas}(u(E))$, we conclude that 
$$
{\rm meas}(u(E)) \le 2C_q m\|\nabla u\|_{L^q(A,\R^2)} \eps^{\frac{q-1}{q}}.
$$
Since $\eps$ is arbitrary, we deduce that ${\rm meas}(u(E))=0$.
\end{proof}

\begin{lemma}
\label{image}
Let $\phi \in W^{1,p}(\R^N) \cap C^0(\R^N)$ with $p \in [1,+\infty]$.
Let $K \subseteq \R^N$ be such that $\phi(K)$ is compact and ${\rm meas}(\phi(K))=0$.
Then there exists $\phi_n \in W^{1,p}(\R^N) \cap C^0(\R^N)$ with
\begin{equation}
\label{strongw1plemma}
\phi_n \to \phi 
\quad\text{strongly in }W^{1,p}(\R^N)\quad\text{if }1 \le p<+\infty,
\end{equation}
\begin{equation}
\label{weakstar}
(\phi_n,\nabla \phi_n) \weakst (\phi,\nabla \phi) 
\quad\text{\wlystar in }L^{\infty}(\R^N,\R^{N+1})\quad\text{if }p=+\infty,
\end{equation}
and such that $\phi_n$ is locally constant on a neighborhood of $K$, i.e., $\nabla \phi_n=0$ a.e. on a neighborhood $A_n$ of $K$.
\end{lemma}

\begin{proof}
By assumption $C:=\phi(K)$ is compact and such that ${\rm meas}(C)=0$.
Let us set 
$$
C_n:=\left\{y \in \R\,:\,{\rm dist}(y,C) \le \frac{1}{n}\right\}
$$
and
$$
T_n(y):=\int_0^y 1_{\R \setminus C_n}(s)\,ds.
$$
Since ${\rm meas}(C_n) \to 0$ as $n \to +\infty$, we get that 
\begin{equation}
\label{tnid}
T_n \to Id \quad\text{pointwise}.
\end{equation}
Moreover, $T_n$ is 1-Lipschitz, $T'_n=0$ a.e. on $C_n$, and 
\begin{equation}
\label{convTnbis}
T_n'(y) \to 1
\quad\text{for a.e. }y \in \R.
\end{equation}
Let us set
\begin{equation}
\label{defconv}
\phi_n:=T_n \circ \phi.
\end{equation}
We have that $\phi_n \in W^{1,p}(\R^N) \cap C^0(\R^N)$, and by the Chain Rule Formula for Sobolev functions (see for instance \cite[Theorem 3.99]{AFP}) we get for a.e. $x \in \R^N$
\begin{equation}
\label{chainrule}
\nabla \phi_n(x)=T'_n(\phi(x)) \nabla \phi(x)
\end{equation}
(recall that $\nabla \phi=0$ a.e. on $\phi^{-1}(C)$ since $C$ has zero measure \cite[Proposition 3.92]{AFP})).
\par
In view of \eqref{chainrule}, $\nabla \phi_n=0$ on $A_n:=\phi^{-1}(C_n)$ which is a neighborhood of $K$. Moreover, by \eqref{tnid} and \eqref{convTnbis}, we have that \eqref{defconv} and \eqref{chainrule} imply that 
$$
\phi_n \to \phi
\quad\text{and}\quad
\nabla \phi_n \to \nabla \phi
\quad\text{a.e. on }\R^N.
$$
Since $|\phi_n| \le |\phi|$ and $|\nabla \phi_n| \le |\nabla \phi|$, we deduce that
\eqref{strongw1plemma} and \eqref{weakstar} hold.
\end{proof}

The following lemma is very close in spirit to \cite[Lemma 3.6]{DMEP}

\begin{lemma}
\label{helmoltz}
Let $\Om \subseteq \R^2$ be open and bounded, and let $q \ge 2$. Let $\psi \in
L^q(\Om,\R^2)$ be such that
$$
\int_\Om \psi \cdot \nabla u\,dx=0
\quad\text{for every }u \in W^{1,2}(\Om).
$$
Then there exists $\phi \in W^{1,q}(\R^2)$ constant on the connected components of $\Om^c$ (in the case $q=2$ constant $c_2$-quasi everywhere) and such that
$$
\nabla \phi=R\psi,
$$
where $R(a,b):=(-b,a)$ denotes a rotation of $90$ degrees counterclockwise. 
\end{lemma}

\begin{proof}
Let us denote by $K_i$, $i\in I$, the connected components of
$\Om^c$, and let $K_0$ be the unbounded one.
\par
Since $\psi \in L^2(\Om,\R^2)$, by Helmholtz decomposition of the space $L^2(\Om,\R^2)$ (see \cite[Theorem 1.1, Chapter III]{Ga}), there exists $\psi_n \in C^{\infty}_c(\Om,\R^2)$ with
$\Div \psi_n=0$ and
\begin{equation}
\label{psintopsi}
\psi_n \to \psi
\qquad\text{strongly in }L^2(\Om,\R^2).
\end{equation}
By setting $\psi_n=0$ outside $\Om$, we can consider $\psi_n$ as defined on the entire $\R^2$. Let us consider $\varphi_n:=R\psi_n$. Since $\R^2$ is simply connected, and $\varphi_n$ has zero-curl, we have that there exists $\phi_n \in C^\infty(\R^2)$ 
such that
$$
\nabla \phi_n=R\psi_n.
$$
In particular $\nabla \phi_n=0$ on a neighborhood $A^n$ of $\Om^c$, so that for every $i \in I$ there exists $c_i^n \in \R$ such that \begin{equation}
\label{phiinconstant}
\phi_n=c_i^n
\quad
\text{on a neighborhood $A_i^n$ of $K_i$.}
\end{equation}
Since $\phi_n$ is well defined up to a constant, we can assume that $\phi_n=0$ on $K_0$. Let $D$ be a disk centered at the origin and such that $\Omb \subseteq D$. By \eqref{psintopsi} we deduce that there exists $\phi \in W^{1,2}_0(D)$ such that 
$$
\phi_n \to \phi
\quad\text{strongly in }W^{1,2}_0(D).
$$
We have that
$$
\nabla \phi=R\psi \in L^q(\Om,\R^2).
$$
We deduce that $\phi \in W^{1,q}_0(D)$, and in particular $\phi \in W^{1,q}(\R^2)$. Since up to a subsequence $\phi_n \to \phi$ $c_2$-q.e., from \eqref{phiinconstant} we deduce that there exists $c_i \in \R$, $i\in I$, such that
\begin{equation}
\label{c2constant}
\phi=c_i
\quad\text{$c_2$-q.e. on }K_i. 
\end{equation}
\par
In the case $q>2$, we have that $\phi$ is H\"older continuous by Sobolev Embedding Theorem. So by Lemma \ref{continuous}, we get that \eqref{c2constant} implies that $\phi$ is constant on  $K_i$, and the proof is concluded. 
\end{proof}

The following theorem contains the density result for the gradients.

\begin{theorem}
\label{mainthml1p}
Let $1 \le p<2$, and let $\Om \in \as_p(\R^2)$ be an admissible domain. Then for every $u \in W^{1,p}(\Om)$ there exists $(u_n)_{n \in \N}$ sequence in $W^{1,2}(\Om)$ such that
$$
\nabla u_n \to \nabla u
\qquad\text{strongly in }L^p(\Om,\R^2).
$$
\end{theorem}

\begin{proof}
Let $K_i$, $i \in I$, be the connected components of $\Om^c$.
Let us consider
\begin{equation}
\label{defH}
H:=\{\nabla v\,:\, v \in W^{1,2}(\Om)\} \subseteq L^p(\Om,\R^2).
\end{equation}
In order to prove the density result, it suffices to check that for every $u \in W^{1,p}(\Om)$ we have
$$
\nabla u \in \overline H,
$$
where $\overline H$ denotes the closure of $H$ in the $L^p$-norm. Since 
$$
\overline H=(H^\perp)^\perp,
$$
where $(\cdot)^\perp$ denotes the orthogonal space in the sense of Banach spaces, we have to check that
\begin{equation}
\label{mainpointbis}
\nabla u \in (H^\perp)^\perp.
\end{equation}
\par
Our strategy to prove \eqref{mainpointbis} is the following. Firstly we characterize the functions $\psi \in H^\perp$, and then we prove that for every $u \in W^{1,p}(\Om)$ we have the orthogonality condition
\begin{equation}
\label{orthrel}
\int_\Om \psi \cdot \nabla u\,dx=0.
\end{equation}

\vskip10pt\noindent
{\bf Step 1: Characterization of $H^\perp$.}
Let $\psi \in H^\perp \subseteq L^{p'}(\Om,\R^2)$, where $p'>2$ is the conjugate exponent of $p$ ($p'=\frac{p}{p-1}$ if $p \in ]1,2[$, $p'=+\infty$ if $p=1$). By definition of $H^\perp$ we have that for every $v \in W^{1,2}(\Om)$
$$
\int_\Om \psi \cdot \nabla v\,dx=0.
$$
By Lemma \ref{helmoltz}, we deduce that there exists $\phi \in W^{1,p'}(\R^2)$ with $\nabla \phi=R\psi$ ($R(a,b):=(-b,a)$ is the rotation of $90$ degrees counterclockwise), 
and such that for every $i \in I$
\begin{equation}
\label{constantonKi}
\phi=c_i \quad\text{on }K_i
\end{equation}
for suitable $c_i \in \R$.

\vskip10pt\noindent
{\bf Step 2: Checking the orthogonality condition.} In order to conclude the proof, it suffices to check that \eqref{orthrel} holds for every $u \in W^{1,p}(\Om)$. By Step 1, we need to check that
\begin{equation}
\label{orthrel2}
\int_\Om \psi \cdot \nabla u\,dx=-\int_\Om R\nabla \phi \cdot \nabla u\,dx=0,
\end{equation}
where $\phi \in W^{1,p'}(\R^2)$ satisfies \eqref{constantonKi} for some $c_i \in \R$, $i \in I$.
\par
Notice that $\phi(\Om^c)$ is compact. Moreover, since $\Om \in \as_p(\R^2)$, there exists a selection $E$ of the connected components $K_i$ of $\Om^c$ such that
$$
\hs^{\frac{p'-2}{p'-1}}(E)=\hs^{2-p}(E)=0
\quad\text{if }1<p<2
$$
and
$$
\hs^1(E)=0
\quad\text{if }p=1.
$$
By \eqref{constantonKi} we get $\phi(\Om^c)=\phi(E)$, and by Lemma \ref{capzero} we have that ${\rm meas}(\phi(\Om^c))={\rm meas}(\phi(E))=0$. Applying Lemma \ref{image}, there exists $\phi_n \in W^{1,p'}(\R^2)$ such that
\begin{equation}
\label{strongw1pbis}
\phi_n \to \phi
\quad\text{strongly in }W^{1,p'}(\R^2)
\quad\text{if }1<p<2,
\end{equation}
\begin{equation}
\label{weakw1pbis}
(\phi_n,\nabla \phi_n) \weakst (\phi,\nabla \phi)
\quad\text{\wlystar in }L^{\infty}(\R^2,\R^3)
\quad\text{if }p=1,
\end{equation}
and
\begin{equation}
\label{gradzero}
\nabla \phi_n=0
\quad\text{on a neighborhood $A_n$ of }\Om^c.
\end{equation}
Notice that $R\nabla \phi_n$ is divergence-free. Up to reducing $A_n$, we can assume that $\Om \setminus \overline{A_n}$ is regular, and that the support of $R\nabla \phi_n$ is contained in $\Om \setminus \overline{A_n}$. Then we have
$$
\int_\Om R\nabla \phi \cdot \nabla u\,dx=\lim_{n \to +\infty} \int_\Om R\nabla \phi_n \cdot \nabla u\,dx=
\lim_{n \to +\infty} \int_{\Om \setminus \overline{A_n}} R\nabla \phi_n \cdot \nabla u\,dx,
$$
and integrating by parts we deduce that 
$$
\lim_{n \to +\infty} \int_{\Om \setminus \overline{A_n}} R\nabla \phi_n \cdot \nabla u\,dx=
\lim_{n \to +\infty} \int_{\Om \setminus \overline{A_n}} \Div(R\nabla \phi_n)\,u\,dx=0,
$$
so that \eqref{orthrel2} is proved, and the proof is concluded.
\end{proof}

\begin{remark}
\label{importantcase}
{\rm
As mentioned in the Introduction, the density result given by Theorem \ref{mainthml1p} (and the similar result for Sobolev spaces Theorem \ref{mainthmW1p}) is useful to establish a link between stability results for linear and nonlinear Neumann problems. Since stability results usually hold under the assumption of a uniform bound on the number of the connected components of the complements of the varying domains (see Section \ref{stability}), the case $\Om^c$ has a finite number of connected components is the relevant one for the applications. 
\par
In this case, the existence of the function $\phi_n$ satisfying conditions \eqref{strongw1pbis} and \eqref{gradzero} in Step 2 of the proof of Theorem \ref{mainthml1p} can be established more directly without using the approximation Lemma \ref{image} as follows (the case $p=1$  is usually not considered in the study of nonlinear Neumann problems in view of a lack of compactness of $W^{1,1}$).
\par
Let $K_0, K_1,\dots,K_m$ be the connected components of $\Om^c$, where $K_0$ is the unbounded one. Let us consider $\xi_0 \in C^\infty(\R^2)$ and
$\xi_i \in C^\infty_c(\R^2)$, $i=1,\dots,m$ such that 
$\xi_0=1$ on a neighborhood of $K_0$, 
$\xi_i=1$ on a neighborhood of $K_i$, and 
$$
{\rm supp }(\xi_h) \cap {\rm supp }(\xi_k)=\emptyset
\qquad\text{for }h \not=k.
$$
By \cite[Theorem 9.1.3]{AH} for every $i=0,1,\dots,m$ we can find $\phi^i_n \in C^\infty(\R^2)$ with
$$
\phi^i_n=c_i
\quad\text{on a neighborhood of }K_i
$$
and such that
$$
\phi^i_n \to \phi
\quad\text{strongly in }W^{1,p'}(\R^2).
$$
Setting
$$
\phi_n:=\left(1-\sum_{i=0}^m \xi_i\right)\phi+\sum_{i=0}^m \xi_i\phi^i_n,
$$
we get that \eqref{strongw1pbis} and \eqref{gradzero} hold.
}
\end{remark}

\begin{remark}
\label{maly}
{\rm
In the proof of Theorem \ref{mainthml1p} we used the assumption that $\Om$ belongs to the class $\as_p(\R^2)$ in order to apply the approximation Lemma \ref{image} and recover the functions $\phi_n$ satisfying \eqref{strongw1pbis}, \eqref{weakw1pbis} and \eqref{gradzero}. Lemma \ref{image} requires that ${\rm meas}(\phi(\Om^c))=0$, and
for $\Om \in \as_p(\R^2)$ every function $\phi \in W^{1,p'}(\R^2)$ constant on the connected components of $\Om^c$ is such that
${\rm meas}(\phi(\Om^c))=0$. In particular this is the case for the functions we need to approximate, that is $\phi \in W^{1,p'}(\R^2)$ such that $R\nabla \phi \in H^\perp$.
\par
We do not know if for a general $\Om$ we can have ${\rm meas}(\phi(\Om^c))=0$ for any $\phi \in W^{1,p'}(\R^2)$ determining an element of $H^\perp$. For such a $\phi$, by Step 1 (and in view of the proof of Lemma \ref{helmoltz}), we have that there exists a sequence of smooth functions $\phi_n$ such that 
\begin{equation}
\label{malygrad}
\nabla \phi_n=0
\quad\text{on a neighborhood of }\Om^c
\end{equation}
and
\begin{equation}
\label{malyw12strong}
\phi_n \to \phi
\quad\text{strongly in }W^{1,2}(\R^2).
\end{equation}
In particular $\phi_n(\Om^c)$ is finite so that ${\rm meas}(\phi_n(\Om^c))=0$. If this always implies that in the limit ${\rm meas}(\phi(\Om^c))=0$, the fact that $\phi$ is energetically more regular than $\phi_n$, i.e., $\phi \in W^{1,p'}(\R^2)$, plays an essential role. 
\par
We can consider indeed the following example which shows a sequence $(\phi_n)_{n \in \N}$ of smooth functions satisfying \eqref{malygrad} and \eqref{malyw12strong} but with $\phi \in W^{1,2}(\R^2) \cap C(\R^2)$ and such that $\phi(\Om^c)$ is the interval $[-1,1]$. Moreover $\Om^c$ can be chosen such that its connected components admit a selection $E$ with
dimension zero, i.e., $\hs^\alpha(E)=0$ for every $\alpha>0$, so that $\Om \in \as_p(\R^2)$. This example heavily relies on a construction proposed by Mal\'y and Martio \cite{MaMa} in connection with the $N$-property of Sobolev transformations.
\par
Let us consider the square $Q:=]-2,2[ \times ]-2,2[$ in $\R^2$, $J:=\{(t,0)\,:\,-1 \le t \le 1\}$, and $\alpha_n \searrow 0$. Since a point has $c_2$-capacity zero, there are functions $u_m \in C^{\infty}(\R^2)$ such that
$$
u_m \to 0
\quad\text{strongly in }W^{1,2}(\R^2)
$$
and such that $0 \le u_m \le 1$, $0 \in {\rm int}\{u_m=1\}$, and $u_m=0$ outside the ball $B(0,1)$. Let $z_1,z_2 \in J$ and $r_0>0$ be such that the balls $B(z_1,r_0)$ and $B(z_2,r_0)$ are disjoint. Let us set
$$
g_m(x):=\frac{1}{2}u_m\left( \frac{x-z_1}{r_0}\right)-\frac{1}{2}u_m\left( \frac{x-z_2}{r_0}\right).
$$
\par
The functions $\phi_n \in C^\infty(\R^2)$ are constructed as follows. Let $\phi_0$ be the constant function equal to $0$. If $n\ge 1$, let us divide the interval $I:=[-1,1]$ in $n$ intervals $I^n_i$ of length $2/n$; we can find points $x^n_i \in J$ and a radius $r_n$ so small that $nr_n^{\alpha_n} \to 0$ and $\phi_{n-1}$ maps $B(x^n_i,r_n)$ to the middle point of $I^n_i$.
Let $B_n:=\bigcup_{i=1}^n B(x^n_i,r_n)$,
$$
h_{m,n}(x):=
\begin{cases}
2^{-n+1}g_m\left( \frac{x-x^n_i}{r_n}\right)      &\text{if $|x-x^n_i| \le r_n$ for some $i$}\\
0      &\text{otherwise,}
\end{cases}
$$
and let $m_n$ be such that 
$$
\|h_{m_n,n}\|_{W^{1,2}(\R^2)} \le 2^{-n}.
$$
We set 
$$
\phi_n:=\phi_{n-1}+h_{m_n,n},
$$
and we denote by $\phi$ the strong limit in $W^{1,2}(\R^2)$ of $(\phi_n)_{n \in \N}$, which is by construction a Cauchy sequence. Notice that $\phi \in W^{1,2}(\R^2) \cap C(\R^2)$, and that the convergence is also uniform.
\par
Let $\Om:=Q \setminus \bigcap_{n \in \N} B_n$.
We have that $\Om^c=Q^c \cup \bigcap_{n \in \N}B_n$. Since $nr_n^{\alpha_n} \to 0$, 
we have that $\hs^\alpha(\cap_{n \in \N} B_n)=0$ for every $\alpha>0$. As a consequence, the connected components of $\Om^c$ admit a selection $E$ such that $\hs^\alpha(E)=0$ for every $\alpha>0$. In particular $\Om \in \as_p(\R^2)$.
\par
By construction we have that $\phi_n$ is constant on a neighborhood of $\Om^c$ but, since $\phi_n \to \phi$ uniformly, it is easy to see that $\phi(\Om^c)=[-1,1]$. Clearly $\phi$ cannot belong to $W^{1,q}(\R^2)$ for some $q>2$, because otherwise its H\"older continuity would imply ${\rm meas}(\phi(\Om^c))=0$.
}
\end{remark}

We are now in a position to prove the main density result of this paper.

\begin{theorem}
\label{mainthmW1p}
Let $1 \le p<2$, and let $\Om \in \as_p(\R^2)$ be an admissible domain. Then $W^{1,2}(\Om)$ is dense in $W^{1,p}(\Om)$.
\end{theorem}

\begin{proof}
The main ingredients in the proof are Theorem \ref{mainthml1p} and a truncation argument. It is not restrictive to assume that $\Om$ is connected, because we can work on each connected component.
\par
Let $u \in W^{1,p}(\Om)$. The density result will be proved if we show that for every $\eps>0$ we can find $(u_n)_{n \in \N}$ sequence in $W^{1,2}(\Om)$ such that
\begin{equation}
\label{limsupineq}
\limsup_{n \to +\infty} \|u_n-u\|_{W^{1,p}(\Om)} \le e_\eps,
\end{equation}
where $e_\eps \to 0$ as $\eps \to 0$.
\par
It is not restrictive to assume that 
$$
u \in W^{1,p}(\Om) \cap L^\infty(\Om).
$$ 
In fact, if $k>0$ and
$$
T_k(u):=\min\{\max\{u,-k\},k\},
$$
we have $T_k(u) \in W^{1,p}(\Om) \cap L^\infty(\Om)$ and
$T_k(u) \to u$ strongly in $W^{1,p}(\Om)$ as $k \to +\infty$. Then if \eqref{limsupineq} holds for $T_k(u)$, by a diagonal argument it also holds for $u$.
\par
Let $A \subset \subset \Om$, $A$ regular and connected, such that
\begin{equation}
\label{defA}
\|\nabla u\|^p_{L^p(\Om \setminus A)}+\|u\|^p_\infty|\Om \setminus A|<\eps.
\end{equation}
By Theorem \ref{mainthml1p}, there exists $v_n \in W^{1,2}(\Om)$ such that
\begin{equation}
\label{gradientapprox}
\nabla v_n \to \nabla u
\quad\text{strongly in }L^p(\Om,\R^2).
\end{equation}
We claim that, up to adding a constant to $v_n$, we can assume that
\begin{equation}
\label{strongsobolev}
v_n \to u \qquad\text{strongly in }W^{1,p}(A).
\end{equation}
Let us set 
\begin{equation}
\label{defun}
u_n:=\min\{\max\{v_n,-\|u\|_\infty\}, \|u\|_\infty\}.
\end{equation}
Notice that $u_n \in W^{1,2}(\Om)$,
\begin{equation}
\label{newstrong}
u_n \to u \quad\text{strongly in }W^{1,p}(A), 
\end{equation}
and that
\begin{equation}
\label{cfrgradients}
|\nabla u_n| \le |\nabla v_n|
\qquad\text{a.e. in }\Om.
\end{equation}
In view of  \eqref{newstrong}, \eqref{cfrgradients}, \eqref{gradientapprox}, \eqref{defun}, and \eqref{defA}, we deduce that
\begin{multline*}
\limsup_{n \to +\infty}\|u_n-u\|^p_{W^{1,p}(\Om)} \le
\limsup_{n \to +\infty}\|u_n-u\|^p_{W^{1,p}(A)}+
\limsup_{n \to +\infty}\|u_n-u\|^p_{W^{1,p}(\Om \setminus \overline A)}\\
=\limsup_{n \to +\infty}\|u_n-u\|^p_{W^{1,p}(\Om \setminus \overline A)} 
\le \limsup_{n \to +\infty}2^{p-1} \left( \int_{\Om \setminus \overline A}|\nabla u_n|^p+|\nabla u|^p+|u_n|^p+|u|^p\,dx \right) \\
\le 2^p(\|\nabla u\|^p_{L^p(\Om \setminus \overline A)}+\|u\|^p_\infty |\Om \setminus \overline A|) \le 2^p\eps
\end{multline*}
so that \eqref{limsupineq} is proved.
\par
In order to complete the proof, let us check that claim \eqref{strongsobolev} holds.
If
$$
c_n:= \frac{1}{|A|}\int_A v_n\,dx,
$$
since $A$ is regular, by Poincar\'e inequality we have
$$
\tilde v_n=v_n-c_n \quad\text{is bounded in } W^{1,p}(A).
$$
Moreover, by the compact embedding of $W^{1,p}(A)$ in $L^p(A)$, there exists $\tilde v \in L^p(A)$ such that up to a subsequence
$$
\tilde v_n \to \tilde v
\quad\text{strongly in }L^p(A).
$$
Since $\nabla \tilde v_n=\nabla v_n$ on $A$, and in view of \eqref{gradientapprox}, we get that $\nabla \tilde v=\nabla u$ in the sense of distributions on $A$. We deduce that $\tilde v \in W^{1,p}(A)$,
$$
\tilde v_n \to \tilde v
\quad\text{strongly in }W^{1,p}(A),
$$
and since $A$ is connected
$$
\tilde v=u+c_A
$$
for some constant $c_A \in \R$. If we set $\hat{v}_n:=v_n-c_n-c_A$
we get
$$
\hat{v}_n \to u
\quad\text{strongly in }W^{1,p}(A),
$$
so that claim \eqref{strongsobolev} is proved.
\end{proof}

\begin{remark}(\bf {The case $\Om$ unbounded})
\label{unbounded}
{\rm
The density of $W^{1,2}(\Om)$ into $W^{1,p}(\Om)$ when $1 \le p<2$ holds also in the case $\Om$ is unbounded but there exists $r_n \to +\infty$ with $\Om_n:=\Om \cap B(0,r_n) \in \as_p(\R^2)$. In fact, if $u \in W^{1,p}(\Om)$, we have $u \in W^{1,p}(\Om_n)$ so that there exists $v^n_k \in W^{1,2}(\Om_n)$ with $v^n_k \to u$ strongly in $W^{1,p}(\Om_n)$ as $k \to +\infty$. If $\chi_n$ is $C^\infty$ function with
$0 \le \chi_n \le 1$, $\chi_n=1$ on $B(0,r_n/2)$ and $\chi_n=0$ outside $B(0,r_n)$, in order to conclude it suffices to choose $u_n:=v^n_{k_n}\chi_n \in W^{1,2}(\Om)$ for $k_n$ sufficiently large.
}
\end{remark}

\section{Applications}
\label{applicsec}

In this section, we give some applications of Theorem \ref{mainthmW1p}
to stability under boundary variations of nonlinear Neumann problems, and to the optimal cutting of a membrane.
\par
Since we work on domains which are not assumed to be regular (for example they may contain cracks), we will use {\it Deny-Lions spaces}. They are defined as follows.
Let $\Om$ be an open subset of $\R^N$, $p \in [1,\infty[$, and $b \in L^{\infty}(\R^N)$ with $b \ge 0$. Let us set
$$
\leb^{1,p}_b(\Om):=\{u \in L^p_{\rm loc}(\Om)\,:\, \nabla u \in L^p(\Om,\R^N), \int_\Om |u|^pb\,dx<+\infty\}.
$$
We say that $u \rs_b v$ if
$$
\int_\Om [|\nabla u-\nabla v|^p+b|u-v|^p]\,dx=0,
$$
and we set
\begin{equation}
\label{defdenylions}
L^{1,p}_b(\Om):=\leb^{1,p}_b(\Om)/ \rs_b
\end{equation}
endowed with the norm $\|u\|:=\|\nabla u\|_{L^p(\Om,\R^N)}+(\int_\Om|u|^pb\,dx)^{1/p}$.
$L^{1,p}_b(\Om)$ is the {\it Deny-Lions space} of exponent $p$ and weight $b$. In the case $b \equiv 0$, it is usually denoted by $L^{1,p}(\Om)$. 
\par
Notice that $W^{1,p}(\Om) \subseteq  L^{1,p}_b(\Om)$. In the case $b \ge c>0$ and $\Om$ is Lipschitz, we have that equality holds, while if $b$ vanishes on subsets with positive measure or $\Om$ is irregular, the inclusion can be strict (see for example \cite[Section 2.7]{Ma2}).
Moreover $W^{1,p}(\Om)$ is always dense in $L^{1,p}_b(\Om)$, as one can check by truncation. As a consequence, in view of Theorem \ref{mainthmW1p}, we have the following density result.

\begin{proposition}
\label{densityweight}
Let $\Om \in \as_p(\R^2)$ be an admissible domain (see Definition \ref{maindef}), and let $b \in L^{\infty}(\R^2)$ such that $b \ge 0$. Then $W^{1,2}(\Om)$ is dense in $L^{1,p}_b(\Om)$ for $1 \le p<2$.
\end{proposition}

Let
\begin{equation}
\label{defolbis}
\Os_l(\R^2):=\{A \subseteq \R^2 \text{ open }:\, \R^2 \setminus A \text{ has at most $l$ connected components}\}.
\end{equation}
For every $u \in L^{1,p}_b(\Om)$, we denote by $\nabla u1_\Om$ and $u1_\Om$ the extension to zero outside $\Om$ of $\nabla u$ and $u$ respectively.
We will use the following proposition due to Bucur and Varchon (see \cite{Buc} and \cite[Theorem 4.1, Remark 5.2]{Buc2}), which is a sort of Mosco limsup condition (see Section \ref{notations}) for the spaces $L^{1,2}$.

\begin{proposition}
\label{bucurmoscolimsup}
Let $\Om_n$ be a sequence in $\Os_l(\R^2)$ converging to $\Om$ in the Hausdorff complementary topology (see Section \ref{notations}) and such that
$$
{\rm meas}(\Om_n \cap \{b>0\}) \to {\rm meas}(\Om \cap \{b>0\}).
$$
Then for every $u \in L^{1,2}_b(\Om)$ there exists $u_n \in L^{1,2}_b(\Om_n)$ such that
$$
\nabla u_n 1_{\Om_n} \to \nabla u1_\Om
\qquad
\text{strongly in }L^2(\R^2,\R^2)
$$
and
$$
u_n 1_{\Om_n} \to u1_\Om
\qquad
\text{strongly in }L^2_b(\R^2),
$$
where $L^2_b(\R^2)$ denotes the $L^2$-space on $\R^2$ with weight $b$.
\end{proposition}

Propositions \ref{densityweight} and \ref{bucurmoscolimsup} will be our main tools in dealing with stability of nonlinear Neumann problems and with the optimal cutting of a membrane.

\subsection{Stability of nonlinear Neumann problems under boundary variations}
\label{stability}

Let $p \in ]1,+\infty[$, and let $A:\R^2 \times \R^2 \to \R^2$ and $B:\R^2 \times \R \to \R$ be two Carath\'eodory functions such that the following conditions hold: there exist $\alpha \in L^{p'}(\R^2)$ ($p':=\frac{p}{p-1}$), $\beta \in L^1(\R^2)$, $0<c_1\le c_2$ such that for almost every $x \in \R^2$ and for every $\xi, \xi_1, \xi_2 \in \R^2$ with $\xi_1 \not=\xi_2$
\begin{equation}
\label{monotone}
(A(x,\xi_1)-A(x,\xi_2))(\xi_1-\xi_2)>0,
\end{equation}
\begin{equation}
\label{growth}
|A(x,\xi)|\le \alpha(x)+c_2|\xi|^{p-1},
\end{equation}
\begin{equation}
\label{below}
A(x,\xi)\cdot \xi \ge \beta(x)+c_1|\xi|^p.
\end{equation}
We assume that $B$ satisfies \eqref{monotone}, \eqref{growth} and \eqref{below} for almost every $x \in \R^2$, and for all $\xi,\xi_1,\xi_2 \in \R$, with $\xi_1 \not= \xi_2$.
\par
Let $\Om$ be a bounded open subset of $\R^2$. We are interested in the stability under boundary variations of $\Om$ of the elliptic equation
\begin{equation}
\label{pb1}
\begin{cases}
-{\rm div}\,A(x,\nabla u)+B(x,u)=0      & \text{in }\Om, \\
A(x,\nabla u)\cdot \nu=0      & \text{on }\partial \Om,
\end{cases}
\end{equation}
where $\nu$ denotes the outer normal to $\partial \Om$.
Since we do not assume any regularity on $A,B$ and on the boundary of $\Om$, we intend \eqref{pb1} in the usual weak sense of Sobolev spaces. More precisely by a solution of problem \eqref{pb1} we mean a function $u_\Om \in W^{1,p}(\Om)$ such that for every test function $\varphi \in W^{1,p}(\Om)$ we have
$$
\int_\Om [A(x,\nabla u_\Om) \nabla \varphi+B(x,u_\Om)\varphi]\,dx=0.
$$
Existence and uniqueness of a solution to \eqref{pb1} follow by well known results on nonlinear elliptic equations with strictly monotone operators (see for instance \cite{Lio}).
\par
Let $(\Om_n)_{n \in \N}$ be a sequence of uniformly bounded open sets in $\R^2$.  
We say that $\Om$ is stable for the Neumann problems \eqref{pb1} along the sequence $(\Om_n)_{n \in \N}$ if
$$
u_{\Om_n}1_{\Om_n} \to u_{\Om}1_{\Om} \quad\text{strongly in }L^p(\R^2)
$$
and
$$
\nabla u_{\Om_n}1_{\Om_n} \to \nabla u_{\Om}1_{\Om} \quad\text{strongly in }L^p(\R^2,\R^2).
$$
\par
Dal Maso, Ebobisse and Ponsiglione \cite[Theorem 2.3]{DMEP} proved that $\Om$ is stable for problem \eqref{pb1} along $(\Om_n)_{n \in \N}$ for every admissible $A$ and $B$ if and only if the space $W^{1,p}(\Om_n)$ converges in the sense of Mosco to $W^{1,p}(\Om)$ (see Section \ref{notations} for a definition).
\par
If $\Om_n$ is in the class $\Os_l(\R^2)$ defined in \eqref{defolbis}, Bucur and Varchon \cite{Buc} proved that if $\Om_n \to \Om$ in the Hausdorff complementary topology (see Section \ref{notations} for a definition) then $W^{1,2}(\Om_n)$ converges in the sense of Mosco to $W^{1,2}(\Om)$ if and only if ${\rm meas}(\Om_n) \to {\rm meas}(\Om)$.
\par
Dal Maso, Ebobisse and Ponsiglione \cite{DMEP} extend this result to the case $1<p<2$ using a technique of {\it nonlinear harmonic conjugates}, and they prove that in the case $p>2$ the result is in general false.
\par
In the following proposition, using Theorem \ref{mainthmW1p} we prove that the result of \cite{DMEP} can be deduced directly by that of \cite{Buc}.

\begin{proposition}
\label{moscoconv}
Let $(\Om_n)_{n \in \N}$ be a sequence of uniformly bounded sets in $\Os_l(\R^2)$ such that $\Om_n$ converges to $\Om$ in the Hausdorff complementary topology, and let $1<p<2$. 
Then $W^{1,p}(\Om_n)$ converges in the sense of Mosco to $W^{1,p}(\Om)$ (and hence problems \eqref{pb1} are stable) if and only if
\begin{equation}
\label{convmeasprop}
{\rm meas}(\Om_n) \to {\rm meas}(\Om).
\end{equation}
\end{proposition}

\begin{proof}
Let us assume that the Mosco convergence holds. Let $\xi \in \R^2$ with $|\xi|=1$.
Let us consider $u \in W^{1,p}(\Om)$ such that $u(x):=\xi \cdot x$. By $(M1)$-condition, there exists $u_n \in W^{1,p}(\Om_n)$ such that
$$
\nabla u_n 1_{\Om_n} \to \nabla u1_\Om
\quad\text{strongly in }L^p(\R^2,\R^2).
$$
Since $|\nabla u_n 1_{\Om_n}-\nabla u1_\Om|=1$ a.e. on $\Om \setminus \Om_n$, we get
\begin{equation}
\label{firstmeas}
\limsup_{n \to +\infty}{\rm meas}(\Om \setminus \Om_n) \le
\lim_{n \to +\infty} \int_{\R^2} |\nabla u_n 1_{\Om_n}-\nabla u1_\Om|^p\,dx=0.
\end{equation}
Let us consider $u_n \in W^{1,p}(\Om_n)$ such that $u_n(x)=\xi \cdot x$. By $(M2)$-condition, up to a subsequence we have that there exists $u \in W^{1,p}(\Om)$ such that
$$
\nabla u_n1_{\Om_n} \weak \nabla u1_\Om
\quad\text{weakly in }L^p(\R^2,\R^2).
$$
Then, if $D$ is a disk containing $\Om_n$ for every $n$ we get
\begin{equation}
\label{firstmeas2}
\xi\lim_{n \to +\infty}{\rm meas}(\Om_n \setminus \Om)=
\lim_{n \to +\infty} \int_{\R^2} \nabla u_n 1_{\Om_n} 1_{D\setminus \Om}\,dx=
\int_{\R^2} \nabla u1_\Om 1_{D \setminus \Om}\,dx=0.
\end{equation}
Combining \eqref{firstmeas} and \eqref{firstmeas2}, we get that \eqref{convmeasprop} holds.
\par
On the contrary, let us assume that \eqref{convmeasprop} holds, and let us prove the Mosco convergence of $W^{1,p}(\Om_n)$ to $W^{1,p}(\Om)$.
\par
Concerning condition $(M2)$, let $u_k \in W^{1,p}(\Om_{n_k})$ be such that
$$
\nabla u_k 1_{\Om_{n_k}} \weak \Phi
\quad\text{weakly in }L^p(\R^2,\R^2)
$$
and
$$
u_k 1_{\Om_{n_k}} \weak \varphi
\quad\text{weakly in }L^p(\R^2)
$$
for some $\Phi \in L^p(\R^2,\R^2)$ and $\varphi \in L^p(\R^2)$. Clearly, since $\Om_n$ converges to $\Om$ in the Hausdorff complementary topology, we have that
$\Phi=\nabla \varphi$ on $\Om$, so that $u:=(\varphi)_{|\Om} \in W^{1,p}(\Om)$. 
In order to conclude that $(M2)$ holds, we have to prove that $\Phi=\nabla u 1_\Om$ and $\varphi=u1_\Om$. Since $\Om_{n_k} \to \Om$ in the Hausdorff complementary topology and ${\rm meas}(\Om_n)  \to {\rm meas}(\Om)$, we have
$$
1_{\Om_{n_k}} \to 1_\Om \quad\text{strongly in }L^1(\R^2).
$$
Then for all $\eta \in C^\infty_c(\R^2,\R^2)$ we deduce
\begin{multline*}
\int_{\R^2} \Phi \cdot \eta\,dx= 
\lim_{k \to +\infty}  \int_{\R^2} [\nabla u_k 1_{\Om_{n_k}}]\cdot \eta\,dx=
 \lim_{k \to +\infty} \int_{\R^2} [\nabla u_k 1_{\Om_{n_k}}]\cdot [\eta1_{\Om_{n_k}}]\,dx\\
 =\int_{\R^2} \Phi \cdot \eta 1_{\Om}\,dx=\int_{\R^2} \Phi 1_{\Om} \cdot \eta\,dx=\int_{\R^2} \nabla u 1_\Om \cdot \eta\,dx
\end{multline*}
so that $\Phi=\nabla u 1_\Om$. Similarly we can prove that $\varphi=u1_\Om$.
\par
Let us prove condition $(M1)$. Since it is sufficient to approximate functions in a dense subset of $W^{1,p}(\Om)$, since $\Os_l(\R^2) \subseteq \as_p(\R^2)$, we can consider in view of Theorem \ref{mainthmW1p} functions $u \in W^{1,2}(\Om)$. Then by 
Proposition \ref{bucurmoscolimsup} (with $b \equiv 1$) there exists $u_n \in W^{1,2}(\Om_n)$ such that
$$
\|u_n1_{\Om_n}-u1_{\Om}\|_{L^2(\R^2)}+
\|\nabla u_n1_{\Om_n}-\nabla u1_{\Om}\|_{L^2(\R^2,\R^2)} \to 0.
$$
Since $W^{1,2}(\Om_n) \subseteq W^{1,p}(\Om_n)$, and since the $L^2$-norm is stronger that $L^p$-norm on bounded domains ($1<p<2$), we deduce that $(M1)$ holds, and the proof is concluded.
\end{proof}

Let us now consider the following nonlinear Neumann problem
\begin{equation}
\label{pbdeg}
\begin{cases}
-{\rm div}\,A(x,\nabla u)+b(x)|u|^{p-2}u=h      & \text{in }\Om, \\
A(x,\nabla u)\cdot \nu=0      & \text{on }\partial \Om,
\end{cases}
\end{equation}
where $A$ is a Carath\'eodory function satisfying conditions \eqref{monotone}, \eqref{growth}, \eqref{below}, $A(x,0)=0$ for a.e. $x \in \R^2$, $b \in L^\infty(\R^2)$ and $b \ge 0$. In order to guarantee the solvability of \eqref{pbdeg}, we assume moreover that $h=bf+g$ with $f,g \in L^{p'}(\R^2)$, ${\rm supp}\,g \subseteq \Om$, and $\int_C g\,dx=0$ for every connected component $C$ of $\Om$. 
We are interested in problem \eqref{pbdeg} because it introduces some degeneracy with respect to problem \eqref{pb1} as $b$ can vanish on regions of $\Om$ with positive measure. 
\par
By a solution of \eqref{pbdeg} we mean a function 
$u_\Om \in L^{1,p}_b(\Om)$ (or more precisely an equivalence class, see \eqref{defdenylions} for a definition) 
such that for every $\varphi \in L^{1,p}_b(\Om)$
$$
\int_\Om [A(x,\nabla u_\Om)\cdot \nabla \varphi
+b(x)|u_\Om|^{p-2}u_\Om\varphi] \,dx=\int_\Om h\varphi\,dx.
$$
Notice that the integrals appearing in the weak formulation of \eqref{pbdeg} are well defined: in particular notice that, if $U$ is a regular open set such that ${\rm supp}\,g \subseteq U \subseteq \overline{U} \subseteq \Om$, and $\bar \varphi$ denotes the average of $\varphi$ on $U$, by H\"older and Poincar\'e inequalities we get
$$
\left|\int_\Om g\varphi\,dx\right|=\left|\int_U g\varphi\,dx\right|=
\left|\int_U g(\varphi-\bar \varphi)\,dx\right| \le C\|g\|_{L^{p'}(\R^2)}\|\nabla \varphi\|_{L^p(\Om)}.
$$
The existence of a solution $u_\Om$ of \eqref{pbdeg} can be established minimizing on $L^{1,p}_b(\Om)$ the functional
$$
F(u):=\int_\Om [A(x,\nabla u)\cdot\nabla u+b(x)|u|^p-hu] \,dx
$$ 
by means of the Direct Method of the Calculus of Variations.
Uniqueness of the solution follows by strict convexity of $F$.
\par
We say that $\Om$ is stable for the Neumann problems \eqref{pbdeg} along the sequence $(\Om_n)_{n \in \N}$ if
$$
\lim_{n \to +\infty}
\int |\nabla u_{\Om_n}1_{\Om_n}-\nabla u_{\Om}1_{\Om}|^p+
b(x)|u_{\Om_n}1_{\Om_n}-u_\Om1_{\Om}|^p \,dx=0.
$$
\par
The stability of Neumann problems \eqref{pbdeg} has been investigated by Bucur and Varchon in \cite{Buc2} in the case $p=2$ and $A(x,\xi)=\xi$ (but it easily generalizes to $A(x,\xi)=a(x)\xi$ with $a$ giving the correct coercivity). 
\par
The main interest in the stability of \eqref{pbdeg} is that, since $b$ is not assumed to be strictly positive, stability is not equivalent to the Mosco convergence of $L^{1,2}_b(\Om_n)$ to 
$L^{1,2}_b(\Om)$ (see \cite[Remark 5.2]{Buc2}).  As a consequence, passing to the nonlinear setting with $1<p<2$, an approach to stability in the line of Dal Maso, Ebobisse and Ponsiglione \cite{DMEP} based on Mosco convergence cannot be directly used in this situation.
\par
Let $(\Om_n)_{n \in \N}$ be a sequence of uniformly bounded open sets in $\R^2$. 
In the case $p=2$, Bucur and Varchon \cite{Buc2} proved that, if $\Om_n \in \Os_l(\R^2)$
and $\Om_n \to \Om$ in the Hausdorff complementary topology, then stability of \eqref{pbdeg} holds if and only if ${\rm meas}(\Om_n \cap \{b>0\}) \to {\rm meas}(\Om \cap \{b>0\})$.
Proposition \ref{densityweight} permits to extend this result to problems \eqref{pbdeg}.

\begin{proposition}
\label{stabilitynonlinear}
Let $(\Om_n)_{n \in \N}$ be a sequence of uniformly bounded open sets in $\Os_l(\R^2)$ converging to $\Om$ in the Hausdorff complementary topology. 
Then $\Om$ is stable along $(\Om_n)_{n \in \N}$ for the Neumann problems \eqref{pbdeg} if and only if
\begin{equation}
\label{convmeas}
{\rm meas}(\Om_n \cap \{b>0\}) \to {\rm meas}(\Om \cap \{b>0\}).
\end{equation}
\end{proposition}

\begin{proof}
Let us assume that stability holds. Then if we choose $f=1$ and $g=0$ so that $h=b$, we deduce that
$$
u_{\Om_n}=1_{\Om_n}
\quad\text{and}\quad
u_{\Om}=1_{\Om}
$$
and that
\begin{equation}
\label{strong1b}
1_{\Om_n} \to 1_{\Om}
\quad\text{strongly in }L^p_b(\R^2),
\end{equation}
where $L^p_b(\R^2)$ denotes the $L^p$-space with weight $b$.
In particular
\begin{equation}
\label{convzeromis}
{\rm meas}([\Om_n \Delta \Om] \cap \{b>0\}) \to 0,
\end{equation}
where $C \Delta D$ denotes the symmetric difference of $C$ and $D$. 
In fact, if \eqref{convzeromis} does not hold, we have that there exists $\chi \in L^{\infty}(\R^2)$ with $\chi 1_{\{b>0\}}=0$, $\chi \not=0$ and
$$
1_{[\Om_n \Delta \Om] \cap \{b>0\}} \weakst \chi
\quad\text{\wlystar in }L^{\infty}(\R^2).
$$
As a consequence it would be
$$
\int_{\R^2} b 1_{[\Om_n \Delta \Om]}\,dx=\int_{\R^2} b1_{[\Om_n \Delta \Om] \cap \{b>0\}}\,dx
\to \int_{\R^2} b\chi\,dx>0
$$
which is against \eqref{strong1b}.
Since
$$
|{\rm meas}(\Om_n \cap \{b>0\})-{\rm meas}(\Om \cap \{b>0\})| \le
{\rm meas}([\Om_n \Delta \Om] \cap \{b>0\}) \to 0,
$$
we deduce that \eqref{convmeas} holds.
\par
Let us assume now that \eqref{convmeas} holds. Let us set $u_n:=u_{\Om_n}$. Choosing $u_n$ as a test in \eqref{pbdeg}
we deduce that $u_n$ is bounded in $L^{1,p}_b(\R^2)$. Up to a subsequence we have that there exist $\Phi \in L^p(\R^2,\R^2)$ and $\varphi \in L^p_b(\R^2)$ such that
\begin{equation}
\label{weakPhi}
\nabla u_n 1_{\Om_n} \weak \Phi
\quad\text{weakly in }L^p(\R^2,\R^2)
\end{equation}
and
\begin{equation}
\label{weaku*}
u_n 1_{\Om_n} \weak \varphi
\quad\text{weakly in }L^p_b(\R^2).
\end{equation}
By the convergence of $\Om_n$ to $\Om$ in the Hausdorff complementary topology, we have that $\Phi=\nabla \varphi$ on $\Om$. In fact let $\Psi \in C^\infty_c(\Om,\R^2)$ and let $U$ be a regular subset of $\Om$ such that ${\rm supp}(\Psi) \subseteq U \subseteq \overline{U} \subseteq \Om$. Then $U \subseteq \Om_n$ for $n$ large and
$$
\int_U \nabla u_n \cdot \Psi\,dx=-\int_U u_n {\rm div }\Psi\,dx.
$$
Let $c_n$ be the average of $u_n$ on $U$. By Poincar\'e inequality and Rellich Compact Embedding of $W^{1,p}(U)$ into $L^p(U)$ we get up to a further subsequence
\begin{equation}
\label{convutilde}
(u_n-c_n) \to \tilde u
\quad\text{strongly in }L^p(U).
\end{equation}
In particular the convergence is strong in $L^p_b(U)$ and
\begin{equation}
\label{weakder}
\int_U \Phi \cdot \Psi\,dx=-\int_U \tilde u\,{\rm div }\Psi\,dx.
\end{equation}
By \eqref{weaku*} and \eqref{convutilde} we deduce that $c_n=u_n-(u_n-c_n)$ converges to some $c \in \R$. We conclude that $\tilde u=\varphi-c$ and by \eqref{weakder} we get
$$
\int_U \Phi \cdot \Psi\,dx=-\int_U \varphi\,{\rm div }\Psi\,dx
$$
which means that $\Phi=\nabla \varphi$ on $\Om$.
\par
Notice moreover that \eqref{convmeas} and \eqref{weaku*} imply that 
\begin{equation}
\label{convintegrale}
\lim_{n \to +\infty}\int_{\R^2} u_n1_{\Om_n}h\,dx=\int_\Om \varphi h\,dx.
\end{equation}
In fact, since $h=bf+g$ and ${\rm supp}(g) \subseteq \Om$, it suffices to check that
$$
\lim_{n \to +\infty}\int_{\R^2} u_n1_{\Om_n}bf\,dx=\int_\Om \varphi bf\,dx.
$$
Since $1_{\Om_n \cap \{b>0\}} \to 1_{\Om \cap  \{b>0\}}$ strongly in $L^1(\R^2)$ we have
$$
\lim_{n \to +\infty}\int_{\R^2} u_n1_{\Om_n}bf\,dx=\lim_{n \to +\infty}\int_{\R^2} u_n1_{\Om_n}bf 1_{\Om_n \cap \{b>0\}}\,dx=\int_{\R^2} \varphi b f1_{\Om \cap \{b>0\}}\,dx.
$$
\par
Let us prove that $\varphi=u_\Om$ on $\Om$. In fact, for every $v \in L^{1,p}_b(\Om)$, by monotonicity we have that
\begin{multline}
\label{214DalMaso}
\int_{\R^2} [A(x,\nabla v1_\Om)\cdot(\nabla v1_\Om-\nabla u_n1_{\Om_n}) 
+B(x,v1_\Om)(v1_\Om-u_n1_{\Om_n})]\,dx \\
\ge 
\int_{\R^2} [A(x,\nabla u_n1_{\Om_n})\cdot(\nabla v1_\Om-\nabla u_n1_{\Om_n})
+B(x,u_n1_{\Om_n})(v1_\Om-u_n1_{\Om_n})]\,dx,
\end{multline}
where
$$
B(x,\xi):=b(x)|\xi|^{p-2}\xi-h(x).
$$
We claim that there exists $v_n \in L^{1,p}_b(\Om_n)$ such that
\begin{equation}
\label{strongvn}
\nabla v_n 1_{\Om_n} \to \nabla v 1_\Om 
\quad\text{strongly in }L^p(\R^2,\R^2)
\end{equation}
and
\begin{equation}
\label{strongvn2}
v_n 1_{\Om_n} \to v 1_\Om 
\quad\text{strongly in }L^p_b(\R^2).
\end{equation}
Notice that
\begin{equation}
\label{convh}
\lim_{n \to +\infty}\int_{\R^2} hv_n 1_{\Om_n}\,dx=\int_{\R^2} hv1_\Om\,dx.
\end{equation}
In fact, because of \eqref{strongvn2}, we have
$$
\lim_{n \to +\infty}\int_{\R^2} bfv_n 1_{\Om_n}\,dx=\int_{\R^2} bfv1_\Om\,dx.
$$
Moreover, if $U$ is regular and with ${\rm supp}(g) \subseteq U \subseteq \overline{U} \subseteq \Om$ we have by Poincar\'e inequality
$$
\left|\int_U g(v_n-v)\,dx\right| \le \|g\|_{L^{p'}(U)} \|\nabla v_n-\nabla v\|_{L^p(U,\R^2)} \to 0.
$$
Since $U \subseteq \Om_n$ for $n$ large enough, we conclude that \eqref{convh} holds.
\par
Using $v_n-u_n$ as a test function in \eqref{pbdeg} we can rewrite the right-hand side of \eqref{214DalMaso} as
\begin{multline}
\label{215DM}
\int_{\R^2} [A(x,\nabla u_n1_{\Om_n})\cdot(\nabla v1_\Om-\nabla u_n1_{\Om_n})
+B(x,u_n1_{\Om_n})(v1_\Om-u_n1_{\Om_n})]\,dx \\
=\int_{\R^2} [A(x,\nabla u_n1_{\Om_n})\cdot(\nabla v1_\Om-\nabla v_n1_{\Om_n})
+B(x,u_n1_{\Om_n})(v1_\Om-v_n1_{\Om_n})]\,dx. 
\end{multline}
Since $A(x,\nabla u_n1_{\Om_n})$ is bounded in $L^{p'}(\R^2,\R^2)$ and
$|u_n|^{p-2}u_n1_{\Om_n}$ is bounded in $L^{p'}_b(\R^2)$, passing to the limit in \eqref{214DalMaso} and in \eqref{215DM}, by claims \eqref{strongvn} and \eqref{strongvn2}, and in view of \eqref{convintegrale} and of the fact that $A(x,0)=0$, we obtain
\begin{multline}
\label{216DM}
\int_{\Om} [A(x,\nabla v)\cdot(\nabla v-\nabla \varphi)
+B(x,v)(v-\varphi)]\,dx \\
\ge 
\lim_{n \to +\infty}
\int_{\R^2} [A(x,\nabla u_n1_{\Om_n})\cdot(\nabla v1_\Om-\nabla v_n1_{\Om_n})
+B(x,u_n1_{\Om_n})(v1_{\Om}-v_n1_{\Om_n})]\,dx=0.
\end{multline}
Taking $v=\varphi\pm \eps z$ in \eqref{216DM}, with $z \in L^{1,p}_b(\Om)$ and $\eps>0$, dividing by $\eps$, and passing to the limit as $\eps \to 0$, we obtain that $\varphi=u_\Om$ in $\Om$. 
\par
Let us now prove that $\varphi=u_\Om 1_{\Om}$, $\Phi=\nabla u_\Om 1_\Om$, and that the convergences in \eqref{weakPhi} and \eqref{weaku*} are indeed strong. Let us take $v:=u_\Om$ in \eqref{215DM}. Since $A(x,0)=0$ we obtain that
\begin{multline*}
\lim_{n \to +\infty}
\int_{\R^2} (A(x,\nabla u_n1_{\Om_n})-A(x,\nabla u_\Om 1_\Om))\cdot(\nabla u_n1_{\Om_n}-\nabla u_\Om1_{\Om})\,dx\\
+\int_{\R^2}b(x)(|u_n|^{p-2}u_n1_{\Om_n}-|u_\Om|^{p-2}u_\Om1_{\Om})
(u_n1_{\Om_n}-u_\Om1_{\Om})]\,dx=0.
\end{multline*}
By monotonicity we get that each integral tends to zero. Now, the strong convergence of
$\nabla u_n1_{\Om_n}$ to $\nabla u_\Om1_{\Om}$ in $L^p(\R^2,\R^2)$ and of
$u_n1_{\Om_n}$ to $u_\Om1_{\Om}$ in $L^p_b(\R^2)$ is a consequence of \cite[Lemma 2.4]{DMEP}. 
\par
In order to conclude the proof, we have to prove our claim on the existence of $v_n \in L^{1,p}_b(\Om_n)$ satisfying \eqref{strongvn} and \eqref{strongvn2}. By Proposition \ref{bucurmoscolimsup}, every function $z \in L^{1,2}_b(\Om)$ is strong limit of a sequence of functions $z_n \in L^{1,2}_b(\Om_n)$ (with the usual extension to zero outside $\Om$ and $\Om_n$ respectively). Since $\Om \in \as_p(\R^2)$, by Proposition \ref{densityweight} $v$ is a strong limit in $L^{1,p}_b(\Om)$ of functions in $W^{1,2}(\Om)$, which in particular are in $L^{1,2}_b(\Om)$. Then the strong approximability for $v$ that we need follows easily using a diagonal argument.
\end{proof}

\subsection{Nonlinear optimal cutting problem}
\label{nloptimalcut}
In this subsection, we apply Theorem \ref{mainthmW1p} to the problem of optimal cutting for a membrane governed by a nonlinear energy.
\par
Let $\Om \subseteq \R^2$ be open, bounded, and with a Lipschitz boundary, and let $x_1,x_2 \in \Om$. Let us set
$$
\ks(\Omb):=\{K \subseteq \Omb\,:\, K \text{ is compact and connected with }x_1,x_2 \in K\}.
$$
Let $f: \R^2 \times \R^2 \to [0,+\infty]$ be a Carath\'eodory function such that $f(x,0)=0$ and satisfying the following growth estimate
\begin{equation}
\label{growthf}
\alpha|\xi|^p \le f(x,\xi) \le \alpha (|\xi|^p+1),
\end{equation}
where $\alpha>0$ and $p \in ]1,+\infty[$.
\par
Let $g \in W^{1,p}(\R^2)$. For every $K \in \ks(\Omb)$ let us set
\begin{equation}
\label{energyfunct}
\Es(K):= \inf \left\{ \int_{\Om \setminus K} f(x,\nabla u)\,dx\,:\, u \in L^{1,p}(\Om \setminus K), u=g \text{ on }\partial \Om \setminus K \right\},
\end{equation}
where $L^{1,p}(\Om \setminus K)$ is the Deny-Lions space defined in \eqref{defdenylions} with $b \equiv 0$. Notice that it is natural to consider displacements in $L^{1,p}(\Om \setminus K)$ because the energy involves only $\nabla u$, and so we cannot expect to control the $L^p$-norm of $u$. Moreover, notice that
the boundary condition on $\partial \Om \setminus K$ is well defined since $\Om$ is Lipschitz
and $u \in W^{1,p}(\Om \cap B_r(x))$ for every $x \in \partial \Om \setminus K$ and $r$ such that $B_r(x) \cap K=\emptyset$.
\par
We have that $\Es(K)$ can be rewritten as
\begin{equation}
\label{newenergy}
\Es(K)=\inf \left\{ \int_{\Om \setminus K} f(x,\nabla u)\,dx\,:\, u \in W^{1,p}(\Om \setminus K), u=g \text{ on }\partial \Om \setminus K \right\}.
\end{equation}
This is due to the density of $W^{1,p}(\Om \setminus K)$ in $L^{1,p}(\Om \setminus K)$, but a little care should be paid for the boundary condition. In particular, denoting by $T_M$ the truncation operator $T_Mu:= \min\{\max\{u,-M\},M\}$, we have 
for every $u \in L^{1,p}(\Om \setminus K)$
$$
\int_{\Om \setminus K} f(x,\nabla u)\,dx=\lim_{M \to +\infty}
\int_{\Om \setminus K} f(x,\nabla T_Mu-\nabla T_M g+\nabla g)\,dx
$$
so that \eqref{newenergy} holds. 
\par
The optimal cutting problem consists in finding the "cut" $K \in \ks(\Omb)$ which maximizes $\Es$ among all admissible cuts, i.e., to solve the problem
\begin{equation}
\label{optimalcutprob}
\max_{K \in \ks(\Omb)}\Es(K).
\end{equation}
The existence of an optimal cutting has been established by Bucur, Buttazzo and Varchon in \cite{BBV} in the case $p=2$ and with a quadratic energy density $f(x,\xi)=A\xi\cdot\xi$. In view of Theorem \ref{mainthmW1p}, the existence of an optimal cut can be proved also in the nonlinear case $1<p<2$. The following result holds.

\begin{proposition}
\label{optimalcutprop}
Let $1<p<2$. Then problem \eqref{optimalcutprob} has a solution.
\end{proposition}

In order to prove Proposition \ref{optimalcutprop} we need the following lemma which is based on Theorem \ref{mainthmW1p} and on the first condition of Mosco convergence for Sobolev spaces $W^{1,2}$ proved in \cite{Buc2}.

\begin{lemma}
\label{densitycut}
Let $1<p<2$, and let $(K_n)_{n \in \N}$ be a sequence in $\ks(\Omb)$ converging in the Hausdorff metric to $K$. Then for every $u \in W^{1,p}(\Om \setminus K)$ with $u=g$ on $\partial \Om \setminus K$ there exists $u_n \in W^{1,p}(\Om \setminus K_n)$ with $u_n=g$ on $\partial \Om \setminus K_n$ such that
$$
\nabla u_n1_{\Om \setminus K_n} \to \nabla u1_{\Om \setminus K}
\quad \text{strongly in }L^p(\R^2,\R^2).
$$
\end{lemma}

\begin{proof}
Let $B$ be an open ball containing $\Omb$. Let us consider
$$
\tilde u:=
\begin{cases}
u      &\text{in }\Om \setminus K\\
g      &\text{in }B \setminus \Omb.
\end{cases}
$$
Notice that $\tilde u \in W^{1,p}(B \setminus K)$. By applying Theorem \ref{mainthmW1p}
to $(B \setminus K) \in \as_p(\R^2)$ we have that there exists $\tilde v_h \in W^{1,2}(B \setminus K)$ such that
$$
\tilde v_h \to \tilde u
\quad\text{strongly in }W^{1,p}(B \setminus K).
$$
By Proposition \ref{bucurmoscolimsup} (with $b \equiv 0$, so that the convergence of the measures is automatically satisfied), for each $h \in \N$ there exists $\tilde v_h^n \in W^{1,2}(B \setminus K_n)$ such that
$$
\nabla \tilde v_h^n1_{B \setminus K_n} \to \nabla \tilde v_h1_{B \setminus K}
\quad\text{strongly in }L^2(\R^2,\R^2).
$$
Using a diagonal argument, and since the $L^2$-convergence is stronger than $L^p$-convergence on bounded domains ($1<p<2$), we can find $\tilde u_n \in W^{1,2}(B \setminus K_n)$ such that
\begin{equation}
\label{nablautildeconvcut}
\nabla \tilde u_n1_{B \setminus K_n} \to \nabla \tilde u1_{B \setminus K}
\quad\text{strongly in }L^p(\R^2,\R^2).
\end{equation}
The functions $\tilde u_n$ do not satisfy a-priori the required boundary condition, so that we have to modify them. We follow here an idea due to Chambolle \cite{Ch}. Up to adding a constant to $\tilde u_n$, we can assume that 
$$
\tilde u_n \to g
\quad
\text{ strongly in }W^{1,p}(B \setminus \Omb).
$$ 
\par
Let us consider $\tilde g_n:=(g-\tilde u_n)_{|B \setminus \Omb}$. We have that
$\tilde g_n \in W^{1,p}(B \setminus \Omb)$ with
$\tilde g_n \to 0$ strongly in $W^{1,p}(B \setminus \Omb)$. Let $E$ denote a linear extension operator from $W^{1,p}(B \setminus \Omb)$ to $W^{1,p}(B)$: such an $E$ exists because $\Om$ has a Lipschitz boundary. Let us set
$$
u_n:=(\tilde u_n+E \tilde g_n)_{|\Om \setminus K_n}.
$$
We have $u_n \in W^{1,p}(\Om \setminus K_n)$ with $u_n=g$ on $\partial \Om \setminus K_n$. Moreover, since $E\tilde g_n \to 0$ strongly in $W^{1,p}(B)$, 
in view of \eqref{nablautildeconvcut}, we deduce that $(u_n)_{n \in \N}$ is the required sequence.
\end{proof}

We are now in a position to prove Proposition \ref{optimalcutprop}.

\begin{proof}[Proof of Proposition \ref{optimalcutprop}]
Let $(K_n)_{n \in \N}$ be a maximizing sequence for the optimal cutting problem, i.e.,
$$
\Es(K_n) \to \sup\{\Es(K)\,:\,K \in \ks(\Omb)\}.
$$
Up to a subsequence, we can assume that $K_n \to K$ in the Hausdorff metric. We have that $K$ is an admissible cut, that is $K \in \ks(\Omb)$. 
\par
By Lemma \ref{densitycut}, for every $u \in W^{1,p}(\Om \setminus K)$ with $u=g$ on $\partial \Om \setminus K$ there exists 
$u_n \in W^{1,p}(\Om \setminus K_n)$ with $u_n=g$ on $\partial \Om \setminus K_n$ such that
$$
\nabla u_n1_{\Om \setminus K_n} \to \nabla u1_{\Om \setminus K}
\quad \text{strongly in }L^p(\R^2,\R^2).
$$
Since $f(x,0)=0$, we deduce
$$
\int_{\Om \setminus K}f(x,\nabla u)\,dx=\lim_{n \to +\infty}
\int_{\Om \setminus K_n}f(x,\nabla u_n)\,dx \ge \lim_{n \to +\infty}\Es(K_n).
$$
Taking the infimum over all admissible $u$ we get
$$
\Es(K) \ge \lim_{n \to +\infty}\Es(K_n)
$$
so that $K$ is an optimal cut, and the proof is concluded.
\end{proof}

\begin{remark}
\label{remcut}
{\rm
In the proof of Proposition \ref{optimalcutprop}, it is not clear if ${\rm meas}(K_n) \to {\rm meas}(K)$ for $n \to +\infty$. As a consequence, the result of Dal Maso, Ebobisse and Ponsiglione \cite{DMEP} cannot be applied to recover the approximability of gradients of functions in $W^{1,p}(\Om \setminus K)$ through gradients of functions in $W^{1,p}(\Om \setminus K_n)$ that we need. It seems essential to use the approximability of the gradients for the relative $W^{1,2}$-spaces established in \cite{Buc2} and the density result given by Theorem \ref{mainthmW1p}.
}
\end{remark}

Let us assume that $f(x,\xi)$ is strictly convex in $\xi$ so that the problem
\begin{equation}
\label{minprobcut}
\min \left\{ \int_{\Om \setminus K} f(x,\nabla u)\,dx\,:\, u \in L^{1,p}(\Om \setminus K), u=g \text{ on }\partial \Om \setminus K \right\}
\end{equation}
admits a unique solution $u_{\Om \setminus K} \in L^{1,p}(\Om \setminus K)$. In particular, in view of \eqref{energyfunct} we have 
$$
\Es(K)=\int_{\Om \setminus K}
f(x,\nabla u_{\Om \setminus K})\,dx.
$$
The associated Euler-Lagrange equation is
\begin{equation}
\label{eulercut}
\begin{cases}
{\rm div}\partial_\xi f(x,\nabla u_{\Om \setminus K})=0      &\text{in }\Om \setminus K\\
\partial_\xi f(x,\nabla u_{\Om \setminus K}) \cdot \nu=0      &\text{on }K\\
u_{\Om \setminus K}=g &\text{on }\partial \Om \setminus K.
\end{cases}
\end{equation}
We deduce that the following stability result for the Neumann-Dirichlet problem \eqref{eulercut} holds.

\begin{proposition}
\label{stabilitycut}
Let $1<p<2$, let $K$ be a solution of \eqref{optimalcutprob}, and let $(K_n)_{n \in \N}$ be a sequence in $\ks(\Omb)$ converging in the Hausdorff metric to $K$. Then we have
that $\Om \setminus K$ is stable for \eqref{eulercut} along the sequence $(\Om \setminus K_n)_{n \in \N}$, that is
\begin{equation}
\label{stabilitycuteq}
\nabla u_{\Om \setminus K_n} 1_{\Om \setminus K_n} \to
\nabla u_{\Om \setminus K} 1_{\Om \setminus K} 
\quad
\text{strongly in }L^p(\R^2,\R^2).
\end{equation}
\end{proposition}

\begin{proof}
Choosing $g$ as an admissible displacement in \eqref{minprobcut}, we get that $\nabla
u_{\Om \setminus K_n} 1_{\Om \setminus K_n}$ is bounded in $L^p(\R^2,\R^2)$ so that up to a subsequence we have
$$
\nabla u_{\Om \setminus K_n} 1_{\Om \setminus K_n} \weak \Phi
\qquad\text{weakly in }L^p(\R^2,\R^2).
$$
Since $K_n \to K$ in the Hausdorff metric, there exists $u \in L^{1,p}(\Om \setminus K)$ with $u=g$ on $\partial \Om \setminus K$ such that $\Phi=\nabla u$ on $\Om \setminus K$. Moreover by lower semicontinuity we have that
\begin{multline}
\label{limsupineqcut}
\int_{\Om \setminus K}f(x,\nabla u)\,dx \le \int_{\R^2} f(x,\Phi)\,dx \le 
\liminf_{n \to +\infty}\int_{\R^2} f(x,\nabla u_{\Om \setminus K_n}1_{\Om \setminus K_n})\,dx\\
=
\liminf_{n \to +\infty}\int_{\Om \setminus K_n}f(x,\nabla u_{\Om \setminus K_n})\,dx.
\end{multline}
Let $v \in L^{1,p}(\Om \setminus K)$ with $v=g$ on $\partial \Om \setminus K$. By Lemma \ref{densitycut}, there exists $v_n \in L^{1,p}(\Om \setminus K_n)$ with $v_n=g$ on $\partial \Om \setminus K_n$ such that
$$
\nabla v_n1_{\Om \setminus K_n} \to
\nabla v1_{\Om \setminus K}
\qquad\text{strongly in }L^p(\R^2,\R^2).
$$
Then by \eqref{limsupineqcut} we get
\begin{multline}
\label{ineqcut2}
\int_{\Om \setminus K}f(x,\nabla v)\,dx=\lim_{n \to +\infty}
\int_{\R^2} f(x,\nabla v_n1_{\Om \setminus K_n})\,dx=
\lim_{n \to +\infty} 
\int_{\Om \setminus K_n}f(x,\nabla v_n)\,dx \\
\ge \liminf_{n \to +\infty} 
\int_{\Om \setminus K_n}f(x,\nabla u_{\Om \setminus K_n})\,dx \ge
\int_{\Om \setminus K}f(x,\nabla u)\,dx
\end{multline}
so that $u=u_{\Om \setminus K}$. Taking $v=u$ in \eqref{ineqcut2} we obtain
$$
\lim_{n \to +\infty}\int_{\Om \setminus K_n}f(x,\nabla u_{\Om \setminus K_n})\,dx=
\int_{\Om \setminus K}f(x,\nabla u_{\Om \setminus K})\,dx.
$$
By \cite{Bre2} we conclude that \eqref{stabilitycuteq} holds, and the proof is concluded.
\end{proof}

\section{Appendix: the density result for the symmetrized gradients}
\label{symsec}
The aim of this appendix is to show how our approach to density explained in Section \ref{secgradients} can be employed to prove the density of the Sobolev space $W^{1,2}$
in the spaces $LD^{1,p}$ of two dimensional elasticity. The case $p=2$ is the really interesting one, and has been proved by Chambolle in \cite{Ch}: using this density, he proves existence for the Cantilever Problem and for the evolution of brittle fractures in the context of planar linearized elasticity. Our approach provides a different proof of Chambolle's result, and covers also the case $1<p<2$.
\par
In order to make the context precise, let $\Om$ be an open subset of $\R^2$. For $1 \le p \le +\infty$, let us set
$$
LD^{1,p}(\Om):=
\left\{ u \in W^{1,p}_{\rm loc}(\Om,\R^2)\,:\, e(u) \in L^p(\Om, \Msym)
\right\},
$$
where $e(u):=(\nabla u+(\nabla u)^T)/2$ denotes the symmetrized gradient of $u$,
and $\Msym$ denotes the space of $2\times 2$ symmetric matrices. 
Clearly $ W^{1,p}(\Om,\R^2) \subseteq  LD^{1,p}(\Om)$. If $\Om$ is Lipschitz, by means of Korn's inequality, it turns out that $LD^{1,p}(\Om)$ coincides with $W^{1,p}(\Om,\R^2)$, while if $\Om$ is irregular, the inclusion can be strict. 
\par
The main result of this section is the following.

\begin{theorem}
\label{mainthmld1p}
Let $1<p \le 2$, and let $\Om \subseteq \R^2$ be a bounded open set  such that $\Om^c$ has a finite number of connected components. Then for every $u \in LD^{1,p}(\Om)$ there exists $u_n \in W^{1,2}(\Om,\R^2)$ such that
$$
e(u_n) \to e(u)
\qquad\text{strongly in }L^p(\Om,\Msym).
$$
\end{theorem}

\begin{proof}
Let $K_i$, $i=0,1,\dots,m$, be the connected components of $\Om^c$, where $K_0$ is the unbounded one. Let us consider the space
$$
H:=\{e(v)\,:\, v \in H^1(\Om,\R^2)\} \subseteq L^p(\Om,\Msym),
$$
where on $\Msym$ we consider the scalar product $A:B:={\rm tr}(AB^T)=\sum_{i,j}a_{ij}b_{ij}$.
\par
In order to prove the theorem, it suffices to check that for every $u \in LD^{1,p}(\Om)$ we have
$$
e(u) \in \overline H,
$$
where the closure is taken in the $L^p$-norm.
\par
We employ a functional analysis argument, namely that $\overline{H}=(H^\perp)^\perp$, where $(\cdot)^\perp$ denotes the orthogonal in the sense of Banach spaces. So our strategy is the following. Firstly we characterize $H^\perp$, and then we check that $e(u)$ is orthogonal to every element of $H^\perp$.

\vskip10pt\noindent
{\bf Step 1: Characterization of $H^\perp$.} 
We claim that
\begin{equation}
\label{Hperpstructsym}
H^\perp:=\left\{\widetilde{\rm Hess}(\varphi)\,:\,\varphi \in W^{2,p'}_0(\R^2), \varphi \text{ is linear on }K_i \text{ for }i=0,1,\dots,m\right\},
\end{equation}
where $\widetilde{\rm Hess}(\varphi)$ is defined as
\begin{equation}
\label{defHessmod}
\widetilde{\rm Hess}(\varphi):=
\left(
\begin{array}{cc}
\partial_2^2\varphi  &-\partial_{12}\varphi \\
-\partial_{12}\varphi &\partial_{1}^2\varphi     
\end{array}
\right).
\end{equation}
By linearity of $\varphi$ on $K_i$ we mean that there exist $c_i \in \R^2$ and $b_i \in \R$ such that (notice that $W^{2,p'}_0(\R^2) \subseteq C^1(\R^2)$)
\begin{equation}
\label{linearityKi}
\varphi(x)=c_i\cdot x+b_i
\qquad\text{for }x \in K_i.
\end{equation}
Since $\varphi \in W^{2,p'}_0(\R^2)$, we clearly have $c_0=0$ and $b_0=0$.
\par
Let us check \eqref{Hperpstructsym}. Let $\Psi \in L^{p'}(\Om,\Msym)$ be an element of $H^\perp$, where $p':=p/(p-1)$ is the conjugate exponent to $p$, with
$$
\Psi=
\left(
\begin{array}{cc}
\psi_1  &\psi_2      \\
\psi_2 &\psi_3     
\end{array}
\right).
$$ 
This means that for every $v \in H^1(\Om,\R^2)$ we have
$$
\int_\Om \Psi\,:\,e(v)\,dx=0.
$$
Choosing $v \in H$ of the form $v:=(v_1,0)$ with $v_1 \in H^1(\Om)$, we deduce that for every $v_1 \in H^1(\Om)$
$$
\int_\Om (\psi_1,\psi_2) \cdot \nabla v_1\,dx=0.
$$
Similarly we deduce that  for every $v_2 \in H^1(\Om)$
$$
\int_\Om (\psi_2,\psi_3) \cdot \nabla v_2\,dx=0.
$$
From Lemma \ref{helmoltz} we conclude that there exist $\phi_1, \phi_2 \in W^{1,p'}(\R^2)$ and $c_i \in \R^2$, $i=0,1,\dots,m$ such that
$$
\nabla \phi_1=R(\psi_1,\psi_2),
\quad
\nabla \phi_2=R(\psi_2,\psi_3),
$$
\begin{equation}
\label{eqci}
(\phi_1,\phi_2)=c_i
\quad\text{on }K_i\text{ for }1<p<2,
\end{equation}
and
\begin{equation}
\label{eqcicap2}
(\phi_1,\phi_2)=c_i
\quad\text{$c_2$-q.e. on }K_i\text{ for }p=2,
\end{equation}
where $R(a,b):=(-b,a)$ denotes a rotation of $90$ degrees counterclockwise.
We can assume $c_0=0$, hence $\phi_1,\phi_2 \in W^{1,p'}_0(\R^2)$.
\par
Let us set
\begin{equation}
\label{defPhi}
\Phi:=(\phi_1,\phi_2) \in W^{1,p'}_0(\R^2,\R^2).
\end{equation}
Let $D$ be a disk centered at the origin and such that $\Omb \subseteq D$. For every $v \in H^1(D,\R^2)$ we have that
$$
\int_D \Phi \cdot \nabla v\,dx=-\int_D (\Div \Phi) \,v\,dx=-\int_D (\partial_1\phi_1+\partial_2\phi_2)v\,dx=-\int_D(-\psi_2+\psi_2)v\,dx=0.
$$
Using again Lemma \ref{helmoltz}, we get that there exists $\varphi \in W^{1,p'}(\R^2)$ with $\varphi=0$ on $\R^2 \setminus D$, and such that
$$
\nabla \varphi=R\Phi=(-\phi_2,\phi_1). 
$$
In view of \eqref{defPhi}, we conclude that $\varphi \in W^{2,p'}_0(\R^2)$. Since $p' \ge 2$, by Sobolev Embedding Theorem we have that $\varphi \in C^1(\R^2)$, so that, by Lemma \ref{continuous}, from \eqref{eqci} and \eqref{eqcicap2} we get (up to replacing $c_i$ with $Rc_i$)
\begin{equation}
\label{eqcibis}
\nabla \varphi=c_i \quad\text{on }K_i.
\end{equation}
\par
By construction we have that $\Psi=\widetilde{\rm Hess}(\varphi)$. In order to complete the proof of claim \eqref{Hperpstructsym}, we need to check \eqref{linearityKi}. Let us consider 
$$
\varphi_i(x):=\varphi(x)-c_i\cdot x.
$$ 
By \eqref{eqcibis}, we clearly have that $\nabla \varphi_i=0$ on $K_i$, i.e., $K_i \subseteq C_i$, where $C_i$ is the set of critical points of $\varphi_i$. By Sard's Lemma we have that 
$$
{\rm meas}(\varphi_i(C_i))=0.
$$
Since $\varphi_i(K_i)$ is connected, and ${\rm meas}(\varphi_i(K_i))=0$, we conclude that $\varphi_i(K_i)=\{b_i\}$ for a suitable $b_i \in \R$, so that \eqref{linearityKi} is proved.

\vskip10pt\noindent
{\bf Step 2: Checking the orthogonality condition.} Let $u \in LD^{1,p}(\Om)$, and let $\Psi \in H^\perp$. We have to check that
$$
\int_\Om \Psi\,:\,e(u)=0.
$$
According to \eqref{Hperpstructsym}, let
$\Psi=\widetilde{\rm Hess}(\varphi)$ with $\varphi \in W^{2,p'}_0(\R^2)$ satisfying \eqref{linearityKi}. Let us consider $\xi_0 \in C^\infty(\R^2)$ and $\xi_i \in C^\infty_c(\R^2)$, $i=1,2,\dots,m$ such that $\xi_0=1$ on a neighborhood of $K_0$, $\xi_i=1$ on a neighborhood of $K_i$, and 
$$
{\rm supp }(\xi_h) \cap {\rm supp }(\xi_k)=\emptyset
\qquad\text{for }h \not=k.
$$
By \cite[Theorem 9.1.3]{AH} we can find $\varphi^i_n \in C^\infty(\R^2)$ with
$$
\varphi^i_n(x)=c_i\cdot x+b_i
\quad\text{on a neighborhood of }K_i,
$$
and such that
$$
\varphi^i_n \to \varphi
\quad\text{strongly in }W^{2,p'}(\R^2).
$$
Let us set
$$
\varphi_n:=\left(1-\sum_{i=0}^m \xi_i\right)\varphi+\sum_{i=0}^m \xi_i\varphi^i_n.
$$
Clearly we have that
\begin{equation}
\label{strongw2pbis}
\varphi_n \to \varphi
\quad\text{strongly in }W^{2,p'}(\R^2),
\end{equation}
and
\begin{equation}
\label{hesszero}
\widetilde{\rm Hess}(\varphi_n)=0
\quad\text{on a neighborhood $A_n$ of }\Om^c.
\end{equation}
We can assume that $\Om \setminus \overline{A_n}$ is regular. Then by means of Korn's inequality we have that  
\begin{equation}
\label{uwpprimo}
u \in W^{1,p}(\Om \setminus \overline{A_n},\R^2).
\end{equation}
By \eqref{strongw2pbis} and \eqref{hesszero} we conclude that
\begin{multline*}
\int_\Om \Psi\,:\,e(u)\,dx=\int_\Om \widetilde{\rm Hess}(\varphi)\,:\,e(u)\,dx=
\lim_{n \to +\infty}\int_\Om \widetilde{\rm Hess}(\varphi_n)\,:\,e(u)\,dx\\
=\lim_{n \to +\infty}\int_{\Om \setminus \overline{A_n}} \widetilde{\rm Hess}(\varphi_n)\,:\,e(u)\,dx. 
\end{multline*}
By \eqref{uwpprimo} and since $\widetilde{\rm Hess}(\varphi_n)$ is symmetric, we deduce that
\begin{equation}
\label{finalequality}
\int_\Om \Psi\,:\,e(u)\,dx=\lim_{n \to +\infty}\int_{\Om \setminus \overline{A_n}} \widetilde{\rm Hess}(\varphi_n)\,:\,\nabla u\,dx. 
\end{equation}
Notice that the rows of $\widetilde{\rm Hess}(\varphi_n)$ are divergence free in $\Om \setminus \overline{A_n}$, and with null trace on $\partial (\Om \setminus \overline{A_n})$. Integrating by parts in \eqref{finalequality}, we get 
$$
\int_\Om \Psi\,:\,e(u)\,dx=0,
$$
so that the proof is concluded.
\end{proof}

\begin{remark}
\label{comparison}
{\rm
In his proof of the density of $W^{1,2}(\Om)$ in $LD^{1,2}(\Om)$, Chambolle \cite{Ch} considers $LD^{1,2}(\Om)$ (up to functions $u$ such that $e(u)=0$) as a Hilbert space with scalar product $(u,v):=\int_\Om e(u):e(v)\,dx$, and proves that
$$
\{e(u)\,:\, u \in W^{1,2}(\Om)\}^{\perp}=0,
$$
where $(\cdot)^\perp$ is the orthogonal in the sense of Hilbert spaces. In this framework, the function $\Psi$ appearing in our Step 1 is of the form $\Psi=e(v)$ for some $v \in LD^{1,2}(\Om)$, and the same analysis implies that $e(v)=\widetilde{{\rm Hess}}(\varphi)$. As a consequence we get $\Delta^2\varphi=0$ ($\varphi$ is usually called the Airy function).
Chambolle uses some PDE and capacity arguments to show that $\varphi=0$
in the case $\Om$ is simply connected, and then proves the general case by reduction to the simply connected one.
\par
In our case, we cannot employ PDE arguments, because we consider $LD^{1,p}(\Om)$ as a natural subspace of $L^p(\Om,\Msym)$, and this seems unavoidable in the case $1<p<2$. As a consequence our function $\varphi$ does not satisfy $\Delta^2\varphi=0$, and we must work out an approximation of $\varphi$ as in Step 2.
}
\end{remark}

\begin{remark}
\label{finitecomponents}
{\rm
In order to follow the arguments of Step 2, it suffices to approximate $\Psi=\widetilde{{\rm Hess}}(\varphi) \in H^\perp$ by $\Psi_n \in L^{p'}(\Om,\Msym)$ which are null on a neighborhood of $\Om^c$ and whose rows are divergence free. This is obtained taking $\Psi_n:=\widetilde{{\rm Hess}}(\varphi_n)$, where $\varphi_n \in W^{2,p'}_0(\R^2)$ is such that
$$
\varphi_n \to \varphi \quad\text{strongly in }W^{2,p'}_0(\R^2),
$$
and with $\varphi_n$ linear on a neighborhood of $\Om^c$. This last constraint cannot be treated using ideas similar to Lemma \ref{image}, so that we used partition of unity (which requires $\Om^c$ with a finite number of connected components) and the approximation result \cite[Theorem 9.1.3]{AH} (which requires $1<p<+\infty$).
}
\end{remark}

\vskip10pt
\noindent {\bf Acknowledgements.}
The authors wish to thank Dorin Bucur for having brought to their attention the problem
of the density of $W^{1,2}$ into $W^{1,p}$ in connection with the stability of nonlinear Neumann problems.

\end{document}